\documentclass[10pt]{article}

\usepackage{moreverb}
\usepackage{libertine}

\usepackage[T1]{fontenc}
\usepackage[english]{babel}
\usepackage[a4paper, left=25mm, right=25mm ]{geometry}

\usepackage{import}
\usepackage{xcolor}
\usepackage[cp1252]{inputenc}

\usepackage[libertine]{newtx}

\usepackage{graphicx}

\usepackage{algorithm}
\usepackage{algorithmicx}
\usepackage{algpseudocode}
\usepackage{mathtools}

\usepackage{booktabs}

\begin{document}

\title{High-Order Numerical Integration on Domains Bounded by Intersecting Level Sets}

\author{\small Lauritz Beck$^{1,2}$, Florian Kummer$^{1,2}$}

\date{ \small
    $^1$ Chair of Fluid Dynamics, Technical University Darmstadt, Germany\\
    $^2$ Graduate School Computational Engineering, Technical University Darmstadt, Germany \\[2ex]
    \today
}

\begin{abstract}
    We present a high-order method
    that provides numerical integration on volumes, 
    surfaces, and lines defined implicitly by two smooth intersecting level sets. 
    To approximate the integrals, 
    the method maps quadrature rules defined on hypercubes to the curved domains of the integrals.
    This enables the numerical integration of a wide range of integrands
    since integration on hypercubes is a well known problem.
    The mappings are constructed by treating the isocontours of the level sets as graphs of height functions.   
    Numerical experiments with smooth integrands indicate a high-order of convergence 
    for transformed Gauss quadrature rules on domains defined 
    by polynomial, rational, and trigonometric level sets.
    We show that the approach we have used can be combined readily with adaptive quadrature methods. 
    Moreover, we apply the approach to numerically integrate on difficult geometries without requiring a low-order fallback method.
\end{abstract}

\maketitle

\section{Introduction}\label{sec1}

In engineering, the natural sciences, or economics,
partial differential equations (PDEs) are prevalently used to model complex processes. 
Because they facilitate numerical solutions with high accuracy, 
high-order methods for solving PDEs have long been a focus of research.
Many high-order numerical methods solving PDEs
are based on discretizing the weak formulations of the PDEs. 
This requires accurate numerical integration 
since the weak formulation involves an integral on the domain of the PDE.
The integrand of the weak formulation belongs to a function space characteristic to the numerical method. 
Common choices for high-order methods are polynomial, rational or trigonometric function spaces \cite{DGreview, IGAreview, gottlieb}.

In multi-phase problems like fluid structure interactions or fluid flows with free surfaces,
separate domains are governed by different PDEs.
A flexible approach to model the geometry of the domains 
is to define the boundaries of the domains implicitly by the zero-isocontours of smooth level sets. 
The level sets often are approximated by functions of the function space of the numerical method.  
For example, this approach is applied in level set methods\cite{olsson,LevelSetMethod} 
or Nitsche type methods\cite{FictitiousDomainMethod,XDG,Reusken}.
Using level sets offers two advantages. First, the boundaries can be moved by simply advecting the level sets. 
Second, a change of topology of the domains is unproblematic for the model.
However, to create a high-order numerical method,
this approach has to be combined with a method 
that provides high-order numerical integration on the domains bounded by level sets.

There exist many multi-phase problems where the boundaries of the domains intersect, 
forming multi-phase contact lines.
A droplet sitting on a soft substrate features a three-phase contact line 
situated at the intersection of the droplet's surface and the substrate's surface, for example.
To model these multi-phase problems numerically, it is sensible to represent each surface by a level set, 
so that the geometry of the contact line is formed by the intersection of their zero-isocontours.
This is particularly suitable,
when the multi-phase problem contains boundary conditions that require accurate geometries, for example surface tension\cite{multiphaseWang}.
Consequently, discretizing the problem requires numerical integration on domains 
that are bounded by intersecting level sets.

Blending multiple smooth level sets
and then applying an integration method for domains bounded by single level sets 
is a straightforward approach to numerically integrate on domains bounded by multiple level sets.
In order to blend multiple level sets, 
they can be combined into a single level set by either multiplying them or by choosing their minimum or maximum.
Although this is a viable approach in combination with robust low-order methods,
the resulting level set generally is not smooth, 
which renders it unfit for most methods of high-order.
Therefore, either a method that can resolve nonsmooth level sets or a method for intersecting level sets is required. 

A variety of high-order methods to integrate on domains bounded by 
level sets have been developed. 
One approach is to extend the domain of integration to a domain for which a quadrature rule is known.
This is achieved by injecting a function into the integral that models the geometry of the domain.
The injected function can be a dirac delta or a Heaviside function\cite{deltaWen1,deltaWen2,deltaWen3},
or it can be defined through a coarea formula\cite{deltaDrescher}.
Respectively relying on error cancelation or a bounded level set, it can be challenging to develop convergent schemes
and importantly, they do not cover nonsmooth geometries or intersecting level sets.

Another approach relies on reconstructing the surfaces to map a known quadrature rule to the domain of the integral.
While some methods refine an initial approximation through perturbation and correction \cite{reconstructionLehrenfeld},
other methods of this approach rely on parametrized 
curves\cite{reconstructionBochkov,reconstructionCheng,reconstructionFries,reconstructionFries2, reconstructionPan}.
Interest in this approach is enhanced by the fact that it overlaps approaches from numerical
methods with explicit surface representations, see\cite{reconstructionEngvall,Antolin} for example.
Some of these methods can resolve kinks, even so reconstructing the implicit surfaces is cumbersome and bounds the accuracy 
of the numerical method by the accuracy of the approximation of the geometry. 

When the geometry of the level set surface can be reconstructed, 
methods allowing unsmooth surfaces can also be derived via integration by parts.
This has been demonstrated for piecewise parameterized boundaries\cite{divGundermann1} or 
trimmed parametric surfaces\cite{divGundermann2}.
In spite of a favorable order of convergence, they also suffer from the disadvantages of the reconstruction approach.

Following a different approach, moment fitting methods construct a quadrature rule by solving a system of equations.
The linear moment fitting methods define quadrature weights through a linear system of equations derived for preset quadrature nodes.
For example, the system of equations can be assembled by applying the divergence theorem \cite{Mueller,divSchwartz}, 
or by utilizing piecewise parameterized boundaries\cite{divSommavaria}. 
Nonlinear methods follow the same approach 
but do not require a preset position of the quadrature nodes\cite{divBui,divHubrich}. 
Although these approaches often result in a comparatively low number of quadrature nodes,
they involve the assembly and solution of a system of equations which is often computationally expensive.
In addition, the approaches are limited to smooth level sets, except for the method relying on reconstruction\cite{divSommavaria}. 
Nevertheless, the reconstruction of the parameterized boundary presents a significant obstacle. 

A unique approach is to recast the zero-isocontour as the graph of an implicit height function \cite{Saye, heightCui}. 
The graph is defined recursively, decomposing the integral into nested, one dimensional integrals. 
Although this offers an efficient high-order integration method for single smooth level sets, 
some geometries require a low-order fallback method
which can reduce efficiency and may deteriorate the convergence order of the method.
If the level set is restricted to polynomials, 
a method without the need for a low-order fallback was presented recently that can resolve geometries of intersecting level sets\cite{SayeNew}. 
However, a quadrature method for domains bounded by general,
smooth intersecting level sets is missing.

In this paper we propose a high-order method
that provides numerical integration on volumes, surfaces, and lines defined implicitly by two intersecting level sets.
Augmenting the approach of Saye\cite{Saye}, the method is based on implicit height functions and restricts the level sets to smooth function spaces only. 
The implicit height functions are used to map hypercubes to the curved domains of the integrals, 
offering a general approach to construct quadrature rules.
Since we conceptually integrate on hypercubes, adaptive quadrature methods can be applied readily.
Moreover, the approach can be applied to integrals defined by single level sets, removing the need for a low-order fallback method.
The method is also available as source code
which is publicly available on github\cite{Beck}.

This paper is organized as follows: After explaining the principle of the method, 
section~\ref{section:gist} introduces its main component, nested mappings.
Section~\ref{section:algorithm} details the two algorithms involved in creating nested mappings 
before utilizing them to derive numerical quadrature rules in section~\ref{section:integration}.
In section~\ref{section:experiments}, numerical experiments are conducted and discussed.
Finally, section~\ref{section:summary} briefly summarizes our results and concludes with further discussions.


\section{The Gist of it} \label{section:gist}

\begin{figure}
    \centering
    \includegraphics{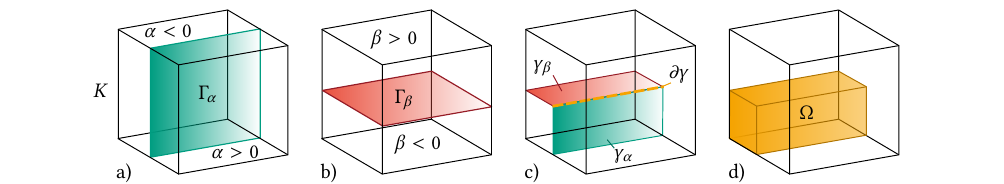}
    \caption{ Implicit definition of the domains of integration. 
    The level sets $\alpha$ and $\beta$ partition hyperrectangle $K$ (a, b),
    defining the surfaces $\Gamma_{\alpha}, \Gamma_{\beta}$, and the boundary lines 
    $\partial \Gamma_{\alpha}, \partial \Gamma_{\beta}$ as subsets of their zero-isocontours (c).
    The region with negative sign of $\alpha$ and $\beta$ specifies the domain of integration $\Omega$ (d). 
    \label{fig:domain}}
\end{figure}

We are interested in numerically integrating a function $g(\boldsymbol{x}),\mathbb{R}^3 \rightarrow \mathbb{R},$ on an implicitly 
defined three-dimensional domain $\Omega \subset \mathbb{R}^3$,
\begin{equation}\label{eq:volume}
    \int_\Omega g(\boldsymbol{x}) \,\text{d}V \approx \sum_i w_i g(\boldsymbol{x}_i),
\end{equation}
with a quadrature rule, consisting of a set of nodes $\boldsymbol{x}_i \in \mathbb{R}^d$ and corresponding weights $w_i \in \mathbb{R}$.
The domain $\Omega$ lies inside a three-dimensional hyperrectangle $K$ 
and is the intersection of the domains $\Omega_{\alpha}$ and $\Omega_{\beta}$,
\begin{equation}
    \Omega = \Omega_{\alpha} \cap \Omega_{\beta}, 
    \qquad \Omega_{\alpha} = \{\boldsymbol{x} \in K | \alpha(\boldsymbol{x}) \leq 0 \}, 
    \quad \Omega_{\beta} = \{ \boldsymbol{x} \in K | \beta(\boldsymbol{x}) \leq 0 \},
\end{equation}
defined by the sign of 
two level set functions $\alpha(\boldsymbol{x}) : K \rightarrow \mathbb{R}$ and 
$\beta(\boldsymbol{x}): K \rightarrow \mathbb{R}$ (see Figure~\ref{fig:domain}).
Further, we are interested in numerically integrating on the two-dimensional 
smooth surfaces $\gamma_{\alpha}$ 
and $\gamma_{\beta}$,
\begin{equation} \label{eq:surface}
    \oint_{\gamma \in \{ \gamma_{\alpha},\gamma_{\beta}\}} g(\boldsymbol{x}) 
        \,\text{d}S,
    \quad
    \gamma_{\alpha} = \Gamma_{\alpha} \cap \Omega_{\beta}, 
    \quad
    \gamma_{\beta} = \Gamma_{\beta} \cap \Omega_{\alpha}, 
\end{equation}
defined by the zero-isocontours of the level set functions,
\begin{equation}
    \Gamma_{\alpha} = \{\boldsymbol{x} \in K | \alpha(\boldsymbol{x}) = 0 \}, 
    \quad
    \Gamma_{\beta} = \{\boldsymbol{x} \in K | \beta(\boldsymbol{x}) = 0 \},
\end{equation}
and on their one-dimensional smooth intersection $\partial \gamma$,
\begin{equation}\label{eq:line}
    \fint_{\partial \gamma} 
        g(\boldsymbol{x}) \,\text{d}L,
    \quad \partial \gamma = \Gamma_{\alpha} \cap \Gamma_{\beta}.
\end{equation}

\subsection{Composition}
The general idea is to simplify the domain of the integral by transforming it into a hypercube.
In case of the volume integral~\eqref{eq:volume}, 
we map a hypercube $C = [-1,1]^3$ to $\Omega$  
with a mapping $T$,
\begin{equation}
    T: C \rightarrow \Omega, \quad T(\tilde{x}, \tilde{y}, \tilde{z}) = (x,y,z), 
\end{equation}
encoding $\Omega$'s geometry in the Jacobian determinant $|\text{det}(D T(\tilde{x}))| = J(T(\tilde{x}))$,
\begin{equation}
    \int_\Omega g(\boldsymbol{x}) \,\text{d}V 
    =  \int_{C} g(T(\boldsymbol{\tilde{x}})) \,|\text{det}(DT(\boldsymbol{\tilde{x}}))| \,\text{d}V
    =  \int_{C} g(T(\boldsymbol{\tilde{x}})) \,J(T(\boldsymbol{\boldsymbol{\tilde{x}}})) \,\text{d}V.
\end{equation}
Since it is a cartesian product of intervals,
a range of quadrature rules suits $C$, e.g. tensorized Gaussian quadrature $\boldsymbol{\tilde{x}}_i, \tilde{w}_i$,
\begin{equation}
    \int_{K} g(T(\boldsymbol{\tilde{x}})) \, J(T(\boldsymbol{\tilde{x}})) \,\text{d}V
    \approx \sum_i \tilde{w}_i g(T(\boldsymbol{\tilde{x}}_i)) J(T(\boldsymbol{\tilde{x}}_i)) .
\end{equation}
After transforming back to $\Omega$,
\begin{equation} 
    \int_\Omega g(\boldsymbol{x}) \,\text{d}V  
    \approx \sum w_i g(\boldsymbol{x}_i), \quad \boldsymbol{x}_i = T(\boldsymbol{\tilde{x}}_i), \quad w_i = J(\boldsymbol{x}_i) \tilde{w}_i.
\end{equation}
we receive a quadrature rule depending only on $T$.

\begin{figure}
    \centering
    \includegraphics{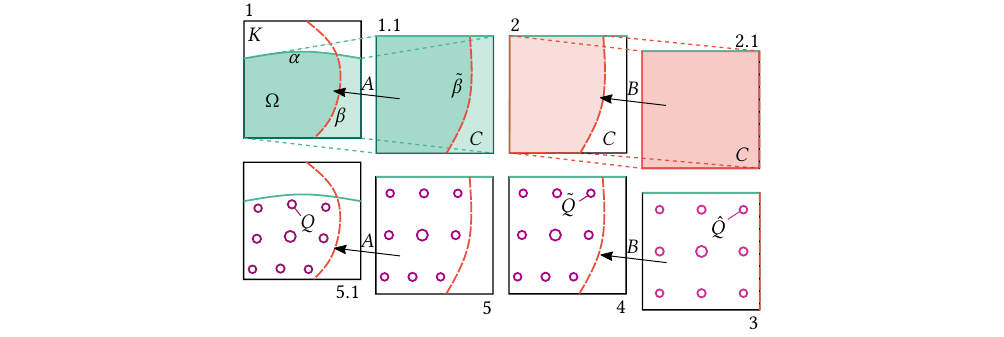}
    \caption{Different steps involved in computing a quadrature rule on domain $\Omega$, colored in dark green (1). 
        After mapping $\Omega$ to a hypercube $C$ in two steps (1-2), 
        the integral may be evaluated via standard tensorized quadrature, e.g. tensorized Gaussian quadrature $\hat{Q}$ (3).
        By embedding $\hat{Q}$ in $\Omega$ with mappings $A$ and $B$ through transforming back (4, 5), 
        the final quadrature rule $Q$ is generated.
        \label{fig:idea}}
\end{figure}

We define $T$ as a composition of mappings for single surfaces, $T = A \circ B$,
utilizing that $\Omega$ 
is bounded by $\Gamma_{\alpha}$ and $\Gamma_{\beta}$.
As illustrated in Figure~\ref{fig:idea}, we proceed in five steps:
\begin{enumerate}
    \item Find a mapping $A: C \rightarrow \Omega_{\alpha}, A(\tilde{x}, \tilde{y}, \tilde{z}) = (x,y,z)$ 
    and map the domain $\Omega_{\alpha}$ defined by level set $\alpha$ to the hypercube $C = [-1,1]^3$, 
        s.t. $\Omega_{\alpha} = A(C)$. 
        Map level set $\beta$ to the hypercube to receive 
        $\tilde{\beta}(\tilde{x}, \tilde{y}, \tilde{z}) 
        = \beta(A(\tilde{x}, \tilde{y}, \tilde{z}))$.
    \item Find a mapping $B: C \rightarrow \Omega_{\tilde{\beta}}, 
        B(\hat{x}, \hat{y}, \hat{z}) = (\tilde{x},\tilde{y},\tilde{z}),$ 
        and map the domain $\Omega_{\tilde{\beta}} = \{\boldsymbol{\tilde{x}} \in C|\tilde{\beta}(\boldsymbol{\tilde{x}}) < 0\}$
        defined by the transformed level set $\tilde{\beta}$ to the hypercube $C = [-1,1]^3$,
        s.t. $\Omega_{\tilde{\beta}} = B(C)$.
    \item Solve the transformed integral numerically on the hypercube $C$
        \begin{equation}
            \int_\Omega g(\boldsymbol{x}) \,\text{d}V 
            = \int_{C} g(A(B(\boldsymbol{\hat{x}}))) \, |\text{det}(D A(B(\boldsymbol{\hat{x}})))| \,\text{d}V 
            \approx \sum_i \hat{w}_i g(T(\boldsymbol{\hat{x}}_i))  J_A(T(\boldsymbol{\hat{x}}_i)) J_B(B(\boldsymbol{\hat{x}}_i))
        \end{equation} 
        using a quadrature rule with nodes $\hat{Q} = \{\boldsymbol{\hat{x}}_1, \boldsymbol{\hat{x}}_2, ...\}$ 
        and weights $\hat{w}_i$.
        The Jacobian determinant of the composition $T$ is separated into its components
        \begin{equation}
            |\text{det}(D A(B(\boldsymbol{\hat{x}})))| = | \text{det}(D A \, D B)| = J_A(T(\boldsymbol{\hat{x}})) J_B(B(\boldsymbol{\hat{x}})).
        \end{equation}
    \item Transform back to receive quadrature nodes 
        $\tilde{Q} = \{\boldsymbol{\tilde{x}}_1, \boldsymbol{\tilde{x}}_2, ...\} = B(\hat{Q}) $ 
        and weights $\tilde{w}_i = \hat{w}_i J_B(\boldsymbol{\tilde{x}}_i)$, 
        \begin{equation}
            \sum_i \hat{w}_i g(T(\boldsymbol{\hat{x}}_i)) J_B(B(\boldsymbol{\hat{x}}_i)) J_A(T(\boldsymbol{\hat{x}}_i))
            = \sum_i \hat{w}_i g(A(\boldsymbol{\tilde{x}}_i)) J_B(\boldsymbol{\tilde{x}}_i) J_A(A(\boldsymbol{\tilde{x}}_i))
            = \sum_i \tilde{w}_i g(A(\boldsymbol{\tilde{x}}_i)) J_A(A(\boldsymbol{\tilde{x}}_i))
            .
        \end{equation}
    \item Transform back quadrature nodes $Q = \{\boldsymbol{x}_1, \boldsymbol{x}_2, ...\} = A(\tilde{Q}) $, 
        yielding the final quadrature rule,
        \begin{equation}
            \int_\Omega g(\boldsymbol{x}) \,\text{d}V 
            \approx \sum_i \tilde{w}_i g(A(\boldsymbol{\tilde{x}}_i)) J_A(A(\boldsymbol{\tilde{x}}_i))
            = \sum_i \tilde{w}_i g(\boldsymbol{x}_i) J_A(\boldsymbol{x}_i)
            = \sum_i w_i g(\boldsymbol{x}_i)
            ,
        \end{equation}
        with nodes $\boldsymbol{x}_i$ and scaled weights $w_i = J_A(\boldsymbol{x}_i) \tilde{w}_i$.
\end{enumerate}

By simply composing mappings, constructing 
a quadrature rule for domains bounded by two level set isocontours is 
reduced to repeatedly mapping a hypercube to domains bounded 
by a single level set isocontour. 
The question that remains is: 
how to find such mappings?

\subsection{Nested Mapping} \label{section_graph}

Representing the backbone of the approach,
let the nested mapping $M: C \rightarrow \Omega_M$,
\begin{equation}
    M(\tilde{x},\tilde{y},\tilde{z}) 
    = 
    \frac{1}{2} 
    \begin{pmatrix}
    (m_x^{\uparrow}  - m_x^{\downarrow}) \tilde{x} + m_x^{\uparrow} + m_x^{\downarrow}\\ 
    (m_y^{\uparrow}(x) - m_y^{\downarrow}(x)) \tilde{y} + m_y^{\uparrow}(x) + m_y^{\downarrow}(x)\\
    (m_z^{\uparrow}(x,y) - m_z^{\downarrow}(x,y)) \tilde{z} + m_z^{\uparrow}(x,y) + m_z^{\downarrow}(x,y)
    \end{pmatrix} 
    =
    \begin{pmatrix}
        x\\
        y\\
        z
    \end{pmatrix}
    ,
\end{equation}
be a composition of height functions,
\begin{equation}
    m_x^{\uparrow}, \quad m_x^{\downarrow}, \quad m_y^{\uparrow} = m_y^{\uparrow}(x), 
    \quad m_y^{\downarrow} = m_y^{\downarrow}(x), \quad m_z^{\uparrow} = m_z^{\uparrow}(x,y), \quad m_z^{\downarrow} = m_z^{\downarrow}(x,y),
\end{equation}
that maps a three-dimensional hypercube $C$ to the codomain $\Omega_M$.
Nesting its components, $M$ is defined recursively:
Evaluating the second entry of $M$ requires its first entry  
and evaluating the third entry of $M$ requires its second entry.

\begin{figure}
    \centering
    \includegraphics{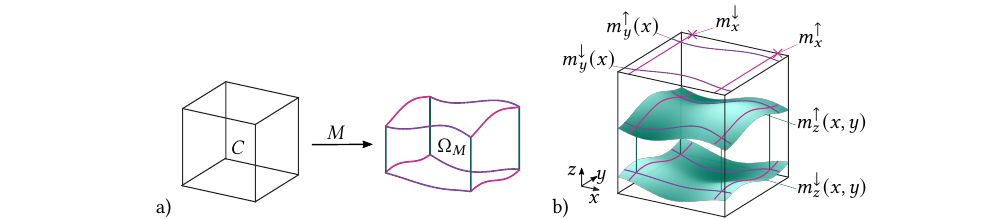}
    \caption{The nested mapping $M$ maps the hypercube $C$ to codomain $\Omega_M$ (a). 
        The codomain is defined by a set of height functions $m_x^{\uparrow}$, 
        $m_x^{\downarrow}$, $m_y^{\uparrow}(x)$, $m_y^{\downarrow}(x)$, 
        $m_z^{\uparrow}(x,y)$, $m_z^{\downarrow}(x,y)$,  
        confining $\Omega_M$ along each coordinate axis (b).
        \label{fig:mappingT}}
\end{figure}

Its height functions do not intersect 
\begin{align}
    m_x^{\downarrow} < m_x^{\uparrow} , 
    \quad m_y^{\downarrow} < m_y^{\uparrow} ,
    \quad m_z^{\downarrow} < m_z^{\uparrow} ,
\end{align}
and confine $\Omega_M$ along the three coordinate axes as shown in Figure~\ref{fig:mappingT}.
First, the constant height functions $m_x^{\uparrow}$ and, $m_x^{\downarrow}$ confine $\Omega_M$ along the $x$-axis. 
Second, the univariate height functions $m_y^{\downarrow}(x)$ and $m_y^{\uparrow}(x)$ confine $\Omega_M$ along the $y$-axis 
and, third, the bivariate height functions $m_z^{\uparrow}(x,y)$ and $m_z^{\downarrow}(x,y)$ confine $\Omega_M$ along the $z$-axis.

Inheriting the nested mappings's properties, the Jacobian matrix $D M$,
\begin{equation}
    D M 
    =
    \frac{1}{2} 
    \begin{pmatrix}
        m^{\uparrow}_x - m^{\downarrow}_x & 0 & 0\\
        M_{xy} \partial_{\tilde{x}} x & m_y^{\uparrow} - m_y^{\downarrow}  & 0 \\
        M_{xz} \partial_{\tilde{x}} x + M_{yz} \partial_{\tilde{x}} y & 
        M_{yz} \partial_{\tilde{y}} y & m_z^{\uparrow} - m_z^{\downarrow}
    \end{pmatrix}
    = 
    \begin{pmatrix}
        \partial_{\tilde{x}} x & 0 & 0\\
        \partial_{\tilde{x}} y & \partial_{\tilde{y}} y & 0\\
        \partial_{\tilde{x}} z & \partial_{\tilde{y}} z & \partial_{\tilde{z}} z
    \end{pmatrix}
    ,
\end{equation}
is structured recursively. 
Its entries, $M_{xy}, M_{xz}, M_{yz}$,
\begin{gather}
    M_{xy} = \partial_x y(x,\tilde{y}) = (\partial_x m_y^{\uparrow} -\partial_x m_y^{\downarrow}) \tilde{y} 
        + \partial_x m_y^{\uparrow} +\partial_x m_y^{\downarrow} ,  \\
    M_{xz} = \partial_x z(x,y, \tilde{z}) = (\partial_x m_z^{\uparrow} - \partial_x m_z^{\downarrow}) \tilde{z}
        + \partial_x m_z^{\uparrow} + \partial_x m_z^{\downarrow},  \\
    M_{yz} = \partial_y z(x,y, \tilde{z}) = ( \partial_y m_z^{\uparrow} - \partial_y m_z^{\downarrow}) \tilde{z} 
        + \partial_y m_z^{\uparrow} + \partial_y m_z^{\downarrow},
\end{gather}
require the first partial derivatives of the height functions.
However, to evaluate the Jacobian determinant $J_M$,
\begin{equation}
    J_M(\tilde{x}, \tilde{y}, \tilde{z}) = \det D M =  
    \frac{1}{8} (m^{\uparrow}_x - m^{\downarrow}_x)  (m_y^{\uparrow} - m_y^{\downarrow} ) (m_z^{\uparrow} - m_z^{\downarrow}),
\end{equation}
partial derivatives are not required.

To increase the possible shapes of $\Omega_M$, the mapping can be oriented differently by rearranging $M$, 
\begin{equation}
    M' = Q^T M(Q(\tilde{x},\tilde{y},\tilde{z})), \quad Q\in \{0,1\}^{3 \times 3}, \quad \det(Q) = 1,
\end{equation}
with a permutation matrix $Q$. 
By rearranging, $\Omega_M$ can be confined along the $x$-axis by a bivariate surface for example.
Naturally, the mapping can be restricted to two-dimensional sets by dropping component $z$.

\subsection{Construction}
Our approach to construct a mapping, $A: C \rightarrow \Omega_\alpha$, 
is to match the codomain of a nested mapping to $\Omega_\alpha$.
To match $\Omega_\alpha$, we disassemble its surface $\Gamma_{\alpha}$ 
into a set of height functions. 
This is possible if $\Gamma_{\alpha}$ and the intersections of $\Gamma_{\alpha}$ with 
suitable faces, edges and corners of the enclosing hyperrectangle $K$,
are graphs on their associated coordinate planes.

\begin{figure}
    \centering
    \includegraphics{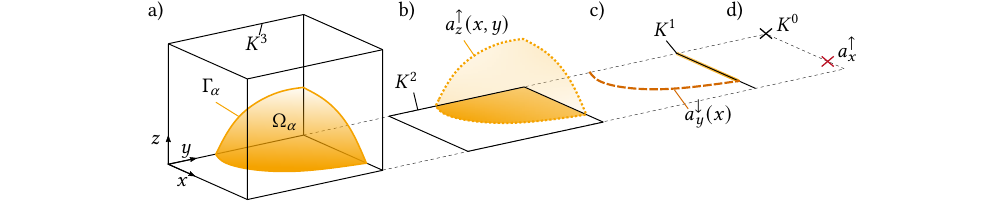}
    \caption{Disassembly of the implicit surface $\Gamma_{\alpha}$ of domain $\Omega_{\alpha}$
        into graphs of the height functions $a_x^{\uparrow}, a_y^{\downarrow}$ and $a_z^{\uparrow}$, 
        from which a nested mapping can be constructed (a):
        Height function $a_z^{\uparrow}(x,y)$ graphs surface $\Gamma_{\alpha} \cap K^3$ over the projection 
            of the surface onto the $xy$-plane (b). 
        Height function $a_y^{\downarrow}$ graphs surface line $\Gamma_{\alpha} \cap K^2$ over the projection
            of the surface line onto of the $xz$-plane (c).
        Height function $a_x^{\uparrow}$ graphs $\Gamma_{\alpha} \cap K^1$ over a point of the $yz$-plane (d).
        \label{fig:mappingA}}
\end{figure}

To illustrate, let $\Omega_\alpha$ be a subset of the cube $K^3 = [-1,1]^3$, 
such that its surface is a graph,
\begin{equation}
    \Gamma_{\alpha} = \left\{ \gamma(\boldsymbol{x}) \,|\, \boldsymbol{x} \in \Pi_{z}(\Gamma_{\alpha}) \right\},
    \quad
    \gamma(\boldsymbol{x}) = 
    \begin{pmatrix}
        x \\
        y \\ 
        a_z^{\uparrow}(x,y)
    \end{pmatrix}
    , 
\end{equation}
on its orthogonal projection $\Pi_{z}(\Gamma_{\alpha})$ onto the $xy$-plane (see figure~\ref{fig:mappingA}).
Hence, by embedding $a^{\uparrow}_z(x,y)$ in $K^3$, the graph's function determines the first height function of the nested mapping.
Forming the graph's domain, 
the orthogonal projection flattens direction $h \in \{x,y,z\}$,
\begin{equation}
    \Pi_h : \mathbb{R}^d \rightarrow  \mathbb{R}^d, \quad \Pi_h(\boldsymbol{x}) = (\mathbb{1} - \boldsymbol{e}_h \otimes \boldsymbol{e}_h) \boldsymbol{x}.
\end{equation}
The projection utilizes height vector $\boldsymbol{e}_h$ which is $1$ at component $h$ and $0$ otherwise.

Further, let the intersection 
of the surface and the cube's lower face $K^2$
be a graph on the orthogonal projection of the intersection onto the $xz$-plane, 
\begin{equation}
    \Gamma_{\alpha} \cap K^2  = \left\{ \gamma(\boldsymbol{x}) \,|\,  \boldsymbol{x} \in \Pi_{y}(\Gamma_{\alpha} \cap K^2) \right\},
    \quad
    \gamma(\boldsymbol{x})  = 
    \begin{pmatrix}
        x \\
        a_y^{\downarrow}(x) \\ 
        z
    \end{pmatrix}
    .
\end{equation}
Repeating the approach, the graph's function $\gamma$ subsequently specifies the second height function of the nested mapping, $a^{\downarrow}_y(x)$.

Finally, the graph of the intersection
of the cube's edge $K^1$ and $\Gamma_{\alpha}$ on its projection onto the $yz$-plane, 
\begin{equation} 
    \Gamma_{\alpha} \cap K^1  = \left\{ \gamma(\boldsymbol{x}) \,|\, \boldsymbol{x} = \Pi_{x}(\Gamma_{\alpha} \cap K^1) \right\},
    \quad
    \gamma(\boldsymbol{x}) = 
    \begin{pmatrix}
        a_x^{\uparrow} \\
        y\\ 
        z
    \end{pmatrix}
    ,
\end{equation}
determines the height function $a^{\uparrow}_x$.

When combined, the height functions define the nested mapping, 
\begin{equation}
    A(\tilde{x},\tilde{y},\tilde{z}) 
    = 
    \frac{1}{2} 
    \begin{pmatrix}
    (a^{\uparrow}_x  + 1) \tilde{x} + a^{\uparrow}_x - 1 \\ 
    (1 - a^{\downarrow}_y(x)) \tilde{y} + 1 + a^{\downarrow}_y(x)\\
    (a^{\uparrow}_z(x,y) + 1) \tilde{z} + a^{\uparrow}_z(x,y) - 1
    \end{pmatrix} 
    =
    \begin{pmatrix}
        x\\
        y\\
        z
    \end{pmatrix}
    .
\end{equation}
where each axis is respectively confined by a height function and a face, edge or point of $K^3$.

\section{Algorithmic Assembly} \label{section:algorithm}

In view of the fact that only mappings for domains bounded by a single level set are required to construct quadrature rules,
this section concentrates on mapping $C$ to $\Omega_{\alpha}$.
Algorithmically, two separate components are involved in finding suitable mappings.
The first component of the algorithm partitions 
$K$ into suitable hyperrectangles, ensuring that $\Gamma_{\alpha}$ is representable by a graph in each.
The second component defines a set of nested mappings 
by tessellating each hyperrectangle with the nested mappings' codomains.

Both components rely on nested bodies; A data structure 
that allows the clear geometric decomposition of a domain bounded by a graph.
A nested body $G = (K_G,H)$ is a strict rooted binary 
tree combining a set of bodies $K_G$ with a set of height directions $H$.
It hierarchically structures the bodies
which are compact sets embedded in $\mathbb{R}^d$, e.g. a square embedded in $\mathbb{R}^3$.
To help distinguish between them, all bodies of equal dimension $i$ are collected in the sets $\mathbb{K}^i$.
Precisely, a nested body $G$ has the following properties:
\begin{enumerate}
    \item Each node of the binary tree is a body $k \in K_G$.
    \item A parent node $k\in \mathbb{K}^i$ with dimension $i$ has two child nodes $k^{\downarrow},k^{\uparrow} \in \mathbb{K}^{i-1}$ with dimension $i-1$.
    \item Each node $k$ of level $l$ of the binary tree is a $d-l$-dimensional body, $k \in \mathbb{K}^{d-l}$.
    \item The nodes of the deepest feasible level $l = d$ consist of zero-dimensional points.
    \item Each level $l < d$ corresponds to a height direction $h_{d-l} \in \{x,y, ...\}$.
    \item The set $H$ gathers the height directions of every level, $H = \{h_{d}, h_{d-1},...\}$.  
\end{enumerate}

\subsection{Graph of the Surface}
\label{subsection_graph}

Algorithm~\ref{alg:graph} recursively creates a nested body $G =(K_G,H)$ consisting of hyperrectangles.
Beginning with a nested body consisting of a single level containing hyperrectangle $K$,
it gradually deepens $G$ by selecting suitable height directions.
Corresponding to the selected height direction, the algorithm adds a new level of hyperrectangles to $G$
while upholding the invariant: For each $k\in K_G $ with associated height direction $h \in H$, 
$\Gamma_{\alpha} \cap k$ is a graph on $\Pi_{h}(\Gamma_{\alpha} \cap k)$.

In each recursion,
Algorithm~\ref{alg:graph} constructs a set of hyperrectangles $P$.
After gathering the lowest level of hyperrectangles of $G$ in $P$,
it removes each hyperrectangle $k$ from $P$ satisfying one of two conditions: 
First, if the hyperrectangle does not intersect $\Gamma_{\alpha}$, $\Gamma_{\alpha} \cap k = \emptyset$.
Second, if the entire hyperrectangle is a subset of $\Gamma_{\alpha}$, $\Gamma_{\alpha} \cap k = k$.
It removes them because the invariant holds unrestrictedly for these hyperrectangles and all contained faces,  
regardless of height direction.
For all remaining hyperrectangles $k\in P$, Algorithm~\ref{alg:graph} chooses a common height direction $h$ and checks 
if the invariant holds.
If it cannot be guaranteed, Algorithm~\ref{alg:graph} subdivides $G$ and restarts or 
designates $G$ for a brute force method.
If the invariant can be guaranteed, it attaches the faces of $P$ as a new level of hyperrectangles to $G$.

Each new level of $G$ consists 
of the lower faces and the upper faces $F$ in direction $h$ of the hyperrectangles $P$,
\begin{equation}
    F =\text{face}^{\downarrow}_{h}(P) \cup \text{face}^{\uparrow}_{h}(P).
\end{equation}
The lower face, $\text{face}_h^{\downarrow}(k)$, 
and upper face, $\text{face}_h^{\uparrow}(k)$,
of a hyperrectangle $k \in \mathbb{K}^i, i > 0,$
\begin{align}
    \text{face}_h^{\downarrow}: \quad \mathbb{K}^i \rightarrow  \mathbb{K}^{i-1}, 
        \quad \text{face}_h^{\downarrow}(k) 
        &= \left\{  \boldsymbol{x} \in k \, | \, 
            \boldsymbol{x} \cdot\boldsymbol{e}_h  \leq \boldsymbol{y} \cdot\boldsymbol{e}_h, 
            \forall \boldsymbol{y} \in k ) \right\}, 
        \\
    \text{face}_h^{\uparrow}: \quad \mathbb{K}^i \rightarrow  \mathbb{K}^{i-1}, 
        \quad \text{face}_h^{\uparrow}(k) 
        &= \left\{  \boldsymbol{x} \in k \, | \, 
            \boldsymbol{x} \cdot\boldsymbol{e}_h  \geq \boldsymbol{y} \cdot\boldsymbol{e}_h, 
            \forall \boldsymbol{y} \in k) \right\}, 
\end{align}
are subsets of the hyperrectangle's boundary.

\begin{figure}
    \centering
    \includegraphics{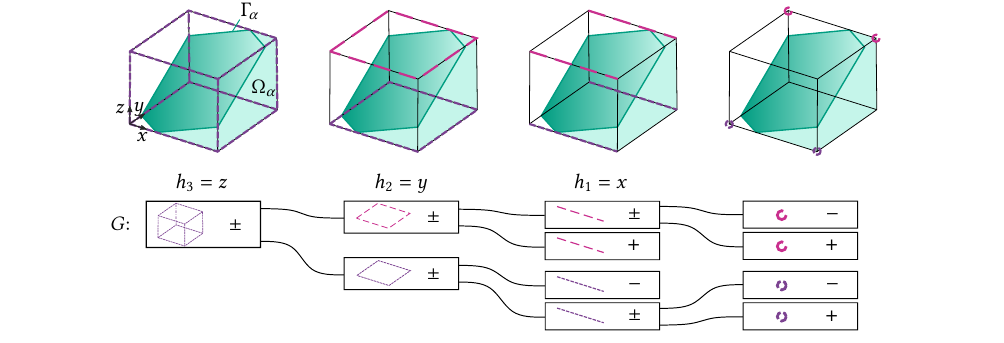}
    \caption{A nested body $G$ created by Algorithm~\ref{alg:graph} for a planar $\Gamma_{\alpha}$, 
        hierarchically linking a hyperrectangle to its faces. 
        Creating a strict rooted binary tree, the faces of each body are attached to each body unless 
        a body has a sign $\in \{-,0,+\}$ imposed by level set $\alpha$. 
        The bodies of each level are successively defined by the level's height direction $h_i$.
        \label{fig:dataTree}}
\end{figure}

As long as $P$ contains hyperrectangles
new levels are attached to $G$ until its deepest level consists of points.
When Algorithm~\ref{alg:graph} has terminated, $\Gamma_{\alpha}$ is a graph whose domain is either 
a hyperrectangle or a subset of a hyperrectangle bounded by graphs.
If the domain is bounded by graphs, 
their domains are in turn either hyperrectangles or subsets of hyperrectangles bounded by graphs, and so on.
Illustrating the involved hyperrectangles, 
Figure~\ref{fig:dataTree} provides an example of a nested body created by Algorithm~\ref{alg:graph} for a planar $\Gamma_{\alpha}$.

\renewcommand{\algorithmicrequire}{\textbf{Input:}}
\renewcommand{\algorithmicensure}{\textbf{Output:}}

\begin{algorithm} 
    \caption{Creates a nested body by gradually adding new levels of faces to it. 
    The surface $\Gamma_{\alpha}$ is a graph on the hyperrectangles of the nested body.}
    \label{alg:graph}
    \begin{algorithmic}[1]
        \Require Nested body $G = (K_G, H)$ and level set $\alpha$. $G$ consists of at least one level of hyperrectangles.
        \Ensure Nested body $G$ consisting of hyperrectangles.
        For each hyperrectangle $k\in K_G$ with height direction $h\in H$ which is not a point,
        $\Gamma_{\alpha} \cap k$ is a graph on $\Pi_{h}(\Gamma_{\alpha} \cap k)$.
        \Procedure{Graph}{$G$, $\alpha$}
            \State $P \gets $ the set of hyperrectangles of the deepest level of $G$
            \State Remove hyperrectangles $\{k \in P | k \cap \Gamma_{\alpha} = k \lor k \cap \Gamma_{\alpha} = \emptyset\}$ from $P$
            \If{$|P| > 0$}
                \State \(h \gets\) height direction of $P$
                \ForAll{$k \in P$}
                    \If{\(\Gamma\sb{\alpha} \cap k\) is not a graph over
                        $\Pi_{h}(\Gamma_{\alpha} \cap k)$}
                        \If{Subdivision count below threshold}
                            \State Split the nested body, \(G \rightarrow G\sp{\prime} , G\sp{\prime\prime}\)
                            \State Restart and call \Call{Graph}{$G\sp{\prime}, \alpha$} 
                            and \Call{Graph}{$G\sp{\prime\prime}, \alpha$}
                        \Else
                            \State Desginate $G$ for a brute force method and \Return 
                        \EndIf
                    \EndIf
                \EndFor
                \State Set the height direction of the deepest level of $G$ to $h$
                \State \(F \gets\) $\text{face}_{h}^{\downarrow}(P) \cup \text{face}_{h}^{\uparrow}(P)$
                \State Add $F$ as new deepest level to \(G\)
                \State \Call{Graph}{$G, \alpha $}
            \EndIf
        \EndProcedure
    \end{algorithmic}   
\end{algorithm}

So far, we have ommitted the inner workings of the algorithm's components
which we will discuss in the following.
To decide if a hyperrectangle should be removed from $P$, 
we assess if $\Gamma_{\alpha} \cap k = \emptyset $ or $\Gamma_{\alpha} \cap k = k$ on the hyperrectangle $k$.
This is equivalent to evaluating the sign $\in \{-,0,+,\pm\}$ of the level set on $k$,
\begin{equation}
    \Gamma_{\alpha} \cap k = \emptyset \Leftrightarrow  \text{sign}_\alpha(k) \in \{-,+\} ,\quad 
    \Gamma_{\alpha} \cap k = k \Leftrightarrow \text{sign}_\alpha(k) = 0.
\end{equation}
The sign,
\begin{equation}
    \text{sign}_\alpha(k)=    
    \begin{cases}
        -, \quad &\text{if } \min\limits_{\boldsymbol{x}\in k} \alpha(\boldsymbol{x}) < 0 
            \land \max\limits_{\boldsymbol{x}\in k} \alpha(\boldsymbol{x}) \leq 0  \\
        0, \quad &\text{if } \min\limits_{\boldsymbol{x}\in k} \alpha(\boldsymbol{x}) = 0 
            \land \max\limits_{\boldsymbol{x}\in k} \alpha(\boldsymbol{x}) = 0\\
        +, \quad &\text{if } \min\limits_{\boldsymbol{x}\in k} \alpha(\boldsymbol{x}) \geq 0 
            \land \max\limits_{\boldsymbol{x}\in k} \alpha(\boldsymbol{x}) > 0 \\
        \pm, \quad &\text{else}
    \end{cases},
\end{equation}
indicates if $\alpha$ can be uniformly bounded on a set $k$, not necessarily a hyperrectangle.

To evaluate the sign numerically, we approximate $\alpha(k)$ with a B\'ezier patch of degree $n$, 
\begin{equation}
    \alpha(k) \approx p(k), \quad p(\boldsymbol{x}) 
    = \sum_i^m q_i b_i(\boldsymbol{u}(\boldsymbol{x})), \quad \boldsymbol{u} : k \rightarrow [-1,1]^{\dim(k)} 
\end{equation}
with $m = (n+1)^{\dim(k)}$ control points $q_i$ and basis functions $b_i$. 
The basis functions are Bernstein polynomials, tensorized to match the dimension of $k$.
The control points are chosen so that $p$ interpolates $\alpha$ at equidistant data points,
spanning a mesh including the boundary of $k$.
Since the B\'ezier patch lies in their convex hull, the control points bound $p(k)$,
\begin{equation}
    \min_i q_i \leq \min_{\boldsymbol{x}\in k} p(\boldsymbol{x}) \leq \max_{\boldsymbol{x}\in k} p(\boldsymbol{x}) \leq \max_i q_i.
\end{equation}
Utilizing the bound, we estimate the sign of $k$,
\begin{equation}
    \text{sign}_\alpha(k) \approx
    \text{sign}_p(k) \approx
    \begin{cases}
        -, \quad &\text{if } \min\limits_i q_i  < -\epsilon \land \max\limits_i q_i \leq \epsilon  \\
        0, \quad &\text{if } \min\limits_i q_i \geq -\epsilon \land \max\limits_i q_i \leq \epsilon\\
        +, \quad &\text{if } \min\limits_i q_i  \geq -\epsilon \land \max\limits_i q_i  > \epsilon \\
        \pm, \quad &\text{else}
    \end{cases}.
\end{equation}
with a precision threshold $\epsilon$, necessitated by the limited precision of machine numbers.
Relying on the convex hull of the approximation, 
the estimate may wrongly assess the body's sign to be undetermined, $\pm$, in some cases.
Fortunately,
this only causes the algorithm to add unnecessary faces to $G$ without compromising the invariant.
In rare cases, the estimate may wrongly determine a sign $\in\{-,0,+\}$ if $\alpha$ lies outside the convex hull of $p$ 
which could compromise correctness. However,
practice has shown that $n = 3$ suffices to reliably assign a sign to a body in a sufficiently resolved numerical grid.

As a way to determine a suitable height direction $h$ for the set of hyperrectangles $P$,
\begin{equation}
    h = \underset{j}{\text{argmax}\,} \boldsymbol{y}_j,
    \quad 
    \boldsymbol{y}  
    = \frac{1}{|P|} \sum_{k \in P} \nabla_k \alpha(\boldsymbol{c}_k),
\end{equation}
the index of the greatest entry of $\boldsymbol{y}$ is selected.
The vector $\boldsymbol{y}$ denotes the average of the projections 
of each gradient $\nabla_k \alpha$ onto each hyperrectangle $k \in P$. 
Each projected gradient is evaluated at the center $\boldsymbol{c}_k$ of its associated hyperrectangle.  

To decide if $\Gamma_\alpha \cap k$ is a graph over $\Pi_h(\Gamma_{\alpha} \cap k)$, 
we bound the gradient $\partial_h \alpha$ of the level set.
It is a graph, 
if the component $h$ of $\partial_h \alpha$ does not change its sign on a hyperrectangle $k$, 
\begin{equation}
    \text{sign}_{\partial_h \alpha}(k) \neq \pm.
\end{equation}
Once again, we interpolate $\partial_h \alpha$ with a B\'ezier patch to evaluate its sign.
Offering a conservative approximation, this estimate may deem surfaces,
which are representable as graphs, unrepresentable and lead to unneeded subdivisions.
However, it does not incorrectly evaluate the surface as representable as a graph which would contradict the invariant. 

Finally, if $\Gamma_\alpha \cap k$ is not a graph over $\Pi_h(\Gamma_{\alpha} \cap k)$, 
the algorithm either subdivides $K$ and restarts or 
employs a brute force method, depending on the number of previous subdivisions.
Since a linear polynomial is always a graph over the faces of a hyperrectangle,
approximating $\alpha$ by a linear polynomial represents a convenient brute force method.
If combined with a sufficiently large number of subdivisions,
high accuracy can be maintained.

\subsection{Tessellation of the Domain}

Algorithm~\ref{alg:subdivide} tesselates the root body $R$ of a nested body $G$ created by Algorithm~\ref{alg:graph},
\begin{equation}
    R = \dot\bigcup V, \quad V = \left\{v| \text{sign}_{\alpha}(v)\in \{-,+\} \right\},
\end{equation}
so that each tile $v \in V$ is the codomain of a nested mapping and is either a subset of $R \cap \Omega_{\alpha}$ or $R \setminus \Omega_{\alpha}$.
In order to create $V$, Algorithm~\ref{alg:subdivide} recursively splits $R$
into suitable subsets which are structured in nested bodies.
Importantly, the subset $V^-$ of $V$,
\begin{equation}
    V^- = \left\{v\in V | \text{sign}_{\alpha}(v)=- \right\},
\end{equation}
is directly required for numerical integration since it tesselates the domain of the integral $R \cap \Omega_{\alpha}$. 

\begin{algorithm} 
    \caption{
        Tessellates a hyperrectangle with the codomains of nested mappings and structures each codomain in a nested body $W$.
        The algorithm creates the nested body by successively assigning signs to the nodes of $W$, 
        deriving the sign of a node from its children or by splitting it along $\Gamma_{\alpha}$.
        }
    \label{alg:subdivide}
    \begin{algorithmic}[1]
        \Require Nested body $W = (K_W, H)$ and level set $\alpha$. 
        Every node $w\in K_W$ with height direction $h$ and children $w^{\uparrow}$ and $w^{\downarrow}$ satisfies the following:
            \begin{enumerate}
                \item If $\text{sign}_{\alpha}(w) = \pm$ , then $\text{sign}_{\alpha}(w^{\downarrow})\neq \pm$ and $\text{sign}_{\alpha}(w^{\uparrow}) \neq \pm .$
                \item $\Gamma_{\alpha} \cap w$ is a graph on $\Pi_{h}(\Gamma_{\alpha} \cap w)$.
                \item $\Pi_h(w) = \Pi_h(w^{\downarrow}) = \Pi_h(w^{\uparrow})$.
            \end{enumerate} 
        \Ensure Nested body $W$. Every node $w \in K_W$ has $\text{sign}(w) \in \{-,0,+\}$.
        \Procedure{Tessellate}{$W$, $\alpha$}
            \ForAll{sets of domains $K^i_W$ of $W$ from second deepest level to level 0}
                \ForAll{$w \in K^i_W$ with sign $\pm$}
                    \State $(w, w^{\downarrow}, w^{\uparrow}) \gets \left(w , \,\text{first child of } w, \,\text{second child of } w \right)$
                    \If{the sign of \(w^{\downarrow}\) is \(+\) and the sign of \(w^{\uparrow}\) is \(-\), or vice versa}
                        \State Split the nested domain \(W \rightarrow W\sp{\prime} , W\sp{\prime\prime}\)
                            along $\Gamma_{\alpha} \cap w$
                            so that \(w = w\sp{\prime} \dot \cup w\sp{\prime\prime}\)
                        \State Update height functions of $W'$ and $W''$
                        \State Sign of \(w\sp{\prime} \gets\) sign of \(w^{\downarrow}\)
                        \State Sign of \(w\sp{\prime\prime} \gets \) sign of \(w^{\uparrow}\)
                        \State Restart and call \Call{Tessellate}{$W\sp{\prime}, \alpha$} 
                            and \Call{Tessellate}{$W\sp{\prime\prime}, \alpha$}
                    \Else
                        \State Sign of \(w \gets\) sign of child of $w$ that is not \(0\)
                    \EndIf
                \EndFor
            \EndFor
        \EndProcedure
    \end{algorithmic} 
\end{algorithm}

When Algorithm~\ref{alg:subdivide} is called for the first time,
it receives a nested body $G$ created by Algorithm~\ref{alg:graph} as input, $W= G$.
In each subsequent iteration, Algorithm~\ref{alg:subdivide} processes the nested body $W$.
It traverses the bodies of $W$ with undetermined sign, $\pm$, in reverse level order starting from the second deepest level, 
assigning a sign $\in \{ -, 0, + \}$ to each.
Visiting each undetermined body $w$ of a level before moving upwards,
it derives the sign of $w$ from the signs of its children, splitting $W$ when required.
When finally the root has been assigned a sign, 
it is the codomain of a nested mapping.

The sign of an undetermined body $w \in K_W$ with height direction $h$ 
and children $w^{\downarrow}$ and $w^{\uparrow}$, 
\begin{equation}
    (w, w^{\downarrow} , w^{\uparrow})_W,
\end{equation}
can be determined from the signs of its children if it fullfills three conditions.
First, both children have a determined sign,
\begin{equation}
    \text{sign}_{\alpha}(w^{\downarrow})\neq \pm, \quad \text{sign}_{\alpha}(w^{\uparrow}) \neq \pm .
\end{equation}
Second, $\Gamma_{\alpha} \cap w$ is a graph over 
$\Pi_{h}(\Gamma_{\alpha} \cap w)$.
Third, the projections in height direction $h$ of $w$ and its two children $w^{\downarrow}$ and $w^{\uparrow}$ coincide,
\begin{equation}
    \Pi_h(w) = \Pi_h(w^{\downarrow}) = \Pi_h(w^{\uparrow}).
\end{equation}

If the three conditions hold, 
there are three different cases specifying $w$'s sign:
First, if the body's children $w^{\downarrow}$ and $w^{\uparrow}$ have the same sign $s\in \{ -, + \}$,
then $\Gamma_{\alpha} \cap w = \emptyset$ and $\text{sign}_{\alpha} (w) = s$.
Second, if $w^{\downarrow}$ has a sign $s\in \{ -, + \}$ and $w^{\uparrow}$ has sign $0$,  
then $\Gamma_{\alpha} \cap w = w^{\uparrow}$  and  $\text{sign}_{\alpha}(w) = s$, or vice versa.
Third, if $w^{\downarrow}$ and $w^{\uparrow}$ have a different sign $\in \{ -, + \}$,
then $\Gamma_{\alpha} \cap w$ lies between the children and the sign of $w$ stays undertermined, $\text{sign}_{\alpha}(w) = \pm$.

\begin{figure}
    \centering
    \includegraphics{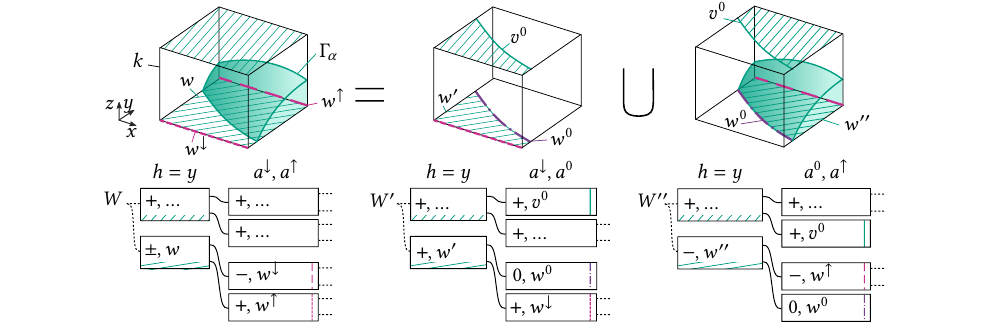}
    \caption{In the first step of splitting the nested body $W$ into $W'$ and $W''$,
        the body $w^0\subseteq \Gamma_{\alpha}$ is inserted into the body $w$. 
        The resulting bodies $w'$ and $w''$ have a sign $\in \{-,+\}$. 
        In the second step, the siblings of $w$ are split by inserting the body $v^0$, 
        which mirrors $w_0$.
        \label{fig:split}}
\end{figure}

To assign a sign in the third case, the algorithm splits $W$ into two nested bodies, $W \rightarrow W', \,W'' $.
It divides $w$ along $\Gamma_{\alpha}$ into two bodies, $w = {w}' \dot \cup {w}''$,
and inserts a new child $w^0 = \Gamma_{\alpha} \cap w$ with $\text{sign}_{\alpha}( w^0) = 0$, 
\begin{equation}
    (w, w^{\downarrow} , w^{\uparrow})_W \rightarrow (w', w^{\downarrow} , w^0)_{W'}, (w'', w^0 , w^{\uparrow})_{W''}.
\end{equation}
The body-children triple of $w$ is split in two,
\begin{equation}
    \Pi_{h}(w') = \Pi_{h}(w^{\downarrow}) = \Pi_{h}(w^0), 
    \quad \Pi_{h}(w'') = \Pi_{h}(w^0) = \Pi_{h}(w^{\uparrow}).
\end{equation}
so that the projections of the new triples align (see Figure~\ref{fig:split}).
After the split, the second case applies to the new bodies ${w}'$ and ${w}''$. 
Hence, $\text{sign}_{\alpha}({w}') = \text{sign}_{\alpha}(w^{\downarrow})$
and $\text{sign}_{\alpha}({w}'') = \text{sign}_{\alpha} (w^{\uparrow})$.

In order to split $W$, not only $w$ has to be divided. To create the other bodies of $W'$ and $W''$, 
the nodes that are part of a level higher than $w$'s are copied and added to $W'$ and $W''$ in a first step. 
In a second step, every remaining node $v$ of $K_W$ is split, $v=v' \dot \cup v''$, 
and inserted into the respective nested body.
Each $v$ is split so that the projections in height direction $h$ of $v$,
\begin{equation}
    \Pi_{h}(v') = \Pi_{h}(v'^{\downarrow}) = \Pi_{h}(v'^{\uparrow}),  \quad 
    \Pi_{h}(v'') = \Pi_{h}(v''^{\downarrow}) = \Pi_{h}(v''^{\uparrow}),
\end{equation}
of $v'$ and its children $v'^{\uparrow}$ and $v'^{\downarrow}$, 
and $v''$ and its children $v''^{\uparrow}$ and $v''^{\downarrow}$ 
align.
In particular, this means that each node $v$ of the level of $w$ receives a new child $v^0$ mirroring $w^0$,
\begin{equation}
    (v, v^{\downarrow} , v^{\uparrow})_W \rightarrow (v', v^{\downarrow} , v^0)_{W'}, (v'', v^0 , v^{\uparrow})_{W''},
\end{equation}
as illustrated by Figure~\ref{fig:split}.

As long as an undetermined body is assigned a sign 
in each recursion and the three conditions are met by every body-children triple of $W$, 
the algorithm will terminate correctly. 
The three conditions are true
when Algorithm~\ref{alg:subdivide} is called for the first time 
because they hold for each nested body created by Algorithm~\ref{alg:graph}.
Consequently, 
all three conditions remain true while the algorithm is iterating and $W$ is not split.
When $W$ is split, it is easy to see that the second and the third condition 
are true for every body-children triple of $W'$ and $W''$; 
The second condition holds
because the bodies of $W'$ and $W''$ are subsets of the bodies of $W$,
and the third condition is enforced by construction.
However, the first condition must be actively maintained by the algorithm.

To assert the first condition, each inserted body $v^0$ is assigned a sign after the split. 
In order to determine the sign, 
it suffices to evaluate $\alpha$ at any point $\boldsymbol{x} \in v^0$. 
This is the only time when explicitly computing a point on the body is required. 
To compute such a point, height functions are used.

Each nested body $W$ has $d$ pairs of height functions,
\begin{equation}
    a_{d-i}^{\downarrow} ,a_{d-i}^{\uparrow}: \mathbb{R}^{d} \rightarrow \mathbb{R}, \quad i = 0, 1, ..., d-1.
\end{equation}
Initially, each height function is constant, 
\begin{equation}
    a_j^{\downarrow}(\boldsymbol{x}) = \min_{\boldsymbol{y} \in R}(e_h \cdot \boldsymbol{y}), 
    \quad a_j^{\uparrow}(\boldsymbol{x}) = \max_{\boldsymbol{y} \in R}(e_h \cdot \boldsymbol{y}), \quad h = h_j,
\end{equation}
and correlates to a face of the hyperrectangle $R$ in the root of $W$.
When a body $w$ of $W$ is split, 
the height functions $a_j^{\downarrow} = a^{\downarrow}, a_j^{\uparrow} = a^{\uparrow}$ 
of $w$'s level are adjusted to match the shape of $W'$ and $W''$
by inserting $a^0 : \mathbb{R}^d \rightarrow \mathbb{R}$,
\begin{equation}
    (a^{\downarrow}, a^{\uparrow})_W \rightarrow  (a^{\downarrow}, a^0)_{W'} \,, (a^0, a^{\uparrow})_{W''}\,.
\end{equation}
The inserted height function $a^0(\boldsymbol{x})$ describes the surface $\Gamma_{\alpha} \cap w$,
\begin{equation}
    \Gamma_{\alpha} \cap w = a^0(w), \quad
    \alpha( \Pi_h(\boldsymbol{x}) + a^0(\boldsymbol{x}) e_h) = 0, \quad  
    a^{\downarrow}(\boldsymbol{x}) \leq a^0(\boldsymbol{x}) \leq a^{\uparrow}(\boldsymbol{x}),
\end{equation}
defined by the zero-isocontour of $\alpha$.
To evaluate the implicit function, the Newton method,
\begin{equation}
    a^0_{n+1} = a^0_n - \frac{\alpha(\boldsymbol{x}_n)}{\partial_{a^0} \alpha(\boldsymbol{x}_n)}  
        = a^0_n - \frac{\alpha(\boldsymbol{x}_n)}{\nabla \alpha(\boldsymbol{x}_n) \cdot e_h}, 
        \quad \boldsymbol{x}_n = \Pi_h(\boldsymbol{x}) + a^0_n(\boldsymbol{x}) e_h,
\end{equation}
is used,
where the subscript $n$ denotes the iteration.
Since $\alpha$ is monotone in height direction $h$ on $w$, the method generally convergences.
Algorithm~\ref{alg:graph} provides monotonicity of $\alpha$ since $w$ is a subset of a face $k$ of $G$, 
\begin{equation}
    \nabla \alpha(x) \cdot e_h \neq 0, \quad \forall x \in w \subseteq k,
\end{equation}
where the gradient's component $h$ can only assume zero on the boundary of $k$.
However, safeguarding is still required, e.g. restricting step-size.

\section{Numerical Integration} \label{section:integration}

\def\mappingsA{A_{\Omega}}
\def\mappingsB{B_{\Omega}}

As described in section~\ref{section:gist}, 
we integrate numerically 
by repeatedly finding nested mappings that map a unit cube 
to domains bounded by single level sets.
To generate a set of nested mappings $\mappingsA$ for $\Omega_{\alpha}$ 
bounded by a single level set,
Algorithm~\ref{alg:graph} first divides $K$ into subsets on which $\Gamma_{\alpha}$ is a graph.
Then Algorithm~\ref{alg:subdivide} 
tessellates $\Omega_{\alpha}$ on each subset,
creating the set of nested bodies $V^-$.
Finally,
the mappings 
$\mappingsA = \{A_0, A_1, ...\}$ are constructed from the set $V^-$ 
utilizing that each nested body $V_i \in V^-$ 
is the codomain of a nested mapping, $A_i : C \rightarrow V_i$. 

The nested mapping $A \in \mappingsA$ can be extracted from the height functions of its nested body constructed by Algorithm~\ref{alg:subdivide}.
Without loss of generality, we assume that the height directions are ordered,
so that $h_1 = x, h_2 = y, h_3 = z$.
Accordingly, the height functions $a_i^{\downarrow}, a_i^{\uparrow}$ returned by Algorithm~\ref{alg:subdivide} 
define the constant, univariate and bivariate height functions,
\begin{gather}
    a_x{\downarrow} = a_1^{\downarrow}(0,0,0), 
    \quad a_y^{\downarrow}(x) = a_2^{\downarrow}(x,0,0), 
    \quad a_z^{\downarrow}(x,y) = a_3^{\downarrow}(x,y,0), \\
    a_x^{\uparrow} = a_1^{\uparrow}(0,0,0), 
    \quad a_y^{\uparrow}(x) = a_2^{\uparrow}(x,0,0), 
    \quad a_z^{\uparrow}(x,y) = a_3^{\uparrow}(x,y,0),
\end{gather}
of $A$, 
\begin{equation}
    A(\tilde{x},\tilde{y},\tilde{z}) 
    = 
    \frac{1}{2} 
    \begin{pmatrix}
    (a_x^{\uparrow}  - a_x^{\downarrow}) \tilde{x} + a_x^{\uparrow} + a_x^{\downarrow}\\ 
    (a_y^{\uparrow}(x) - a_y^{\downarrow}(x)) \tilde{y} + a_y^{\uparrow}(x) + a_y^{\downarrow}(x)\\
    (a_z^{\uparrow}(x,y) - a_z^{\downarrow}(x,y)) \tilde{z} + a_z^{\uparrow}(x,y) + a_z^{\downarrow}(x,y)
    \end{pmatrix} 
    =
    \begin{pmatrix}
        x\\
        y\\
        z
    \end{pmatrix}
    .
\end{equation}

Required in the entries of $A$'s Jacobian matrix, we determine the partial derivatives 
of the height functions describing surface $\Gamma_{\alpha}$ 
by differentiating its implicit definition,
\begin{equation} \label{eq:phi_gamma}
    \frac{\text{d}}{\text{d}x} \left(\alpha(x,y,a_z(x,y)) = 0 \right) 
        \Leftrightarrow \partial_x \alpha
        + \partial_z \alpha \, \partial_x a_z = 0, \quad
    \frac{\text{d}}{\text{d}y} \left(\alpha(x,y,a_z(x,y)) = 0 \right) 
        \Leftrightarrow \partial_y \alpha  + \partial_z \alpha \, \partial_y a_z = 0.
\end{equation}
For brevity,
the height function $a_z$ denotes $a_z^{\uparrow}$ or $a_z^{\downarrow}$.
Rearranging yields the partial derivatives of $a_z$,
\begin{equation}
    \partial_x a_z = - \frac{\partial_x \alpha}{\partial_z \alpha}, \quad
    \partial_y a_z = - \frac{\partial_y \alpha}{\partial_z \alpha},
\end{equation}
which are fractions of the partial derivatives of $\alpha$.
Similarly, differentiating the implicit surface line,
\begin{equation} \label{eq:phi_dgamma}
    \frac{\text{d}}{\text{d}x} \left(\alpha(x,a_y(x),0) = 0 \right)
    \Leftrightarrow \partial_x \alpha + \partial_y \alpha \partial_x a_y  = 0
    \Leftrightarrow \partial_x a_y = - \frac{\partial_x \alpha}{\partial_y \alpha},
\end{equation}
determines the partial derivative $\partial_x a_y$, 
where again $a_y$ denotes $a_y^{\uparrow}$ or $a_y^{\downarrow}$ for brevity.

\subsection{Volume Integrals}

To integrate over domains bounded by a single level set $\alpha$, the volume integral,
\begin{equation}
    \int_{\Omega_{\alpha}} g(\boldsymbol{x}) \,\text{d}V  
    = \sum_{A \in \mappingsA} \int_{C} g(A(\tilde{\boldsymbol{x}})) J_A(\tilde{\boldsymbol{x}})\,\text{d}V 
    \approx \sum_{A \in \mappingsA} \sum_i w_{A,i} g(\boldsymbol{x}_{A,i}),
\end{equation}
is transformed into a sum of volume integrals over hypercube $C$ by the nested mappings of $\mappingsA$.
The resulting volume integrals are evaluated with transformed quadrature nodes $\boldsymbol{x}_{T,i}$ and weights $ w_{T,i}$,
\begin{equation}
    w_{A,i} = J_A(\tilde{x}_i) \tilde{w}_i, \quad \boldsymbol{x}_{A,i} = A(\tilde{\boldsymbol{x}}_i),
\end{equation}
derived from the quadrature nodes $\tilde{\boldsymbol{x}}_i$ and weights $\tilde{w}_i$ of a quadrature rule on $C$.

Applying the same ansatz twice, the volume integral confined by two level sets $\alpha$ and $\beta$,
\begin{equation}
    \int_{\Omega} g(\boldsymbol{x}) \,\text{d}V  
    = \sum_{A \in\mappingsA} \sum_{B \in \mappingsB}\int_{C} g(T(\hat{\boldsymbol{x}})) J_T(\hat{\boldsymbol{x}})\,\text{d}V 
    \approx \sum_{A \in \mappingsA} \sum_{B \in \mappingsB} \sum_i  w_{T,i} g(\boldsymbol{x}_{T,i}),
\end{equation}
is transformed by the composition of mappings $T$,
\begin{equation}
    T: \mathbb{R}^3 \rightarrow \mathbb{R}^3, \quad T = A \circ B .
\end{equation}
The codomains of the nested mappings $\mappingsB =\{B_1, B_2, ...\}$ tesselate $\Omega_{\tilde{\beta}}$ 
defined by $\tilde{\beta}(\tilde{\boldsymbol{x}}) = \beta(A(\tilde{\boldsymbol{x}}))$.
Again, the resulting quadrature nodes and weights,
\begin{equation}
    w_{T,i} = J_T(\hat{\boldsymbol{x}}_i) \hat{w}_i 
    = J_A(B(\hat{\boldsymbol{x}}_i)) J_B(\hat{\boldsymbol{x}}_i) \hat{w}_i,\quad  \boldsymbol{x}_{T,i} 
    = T(\hat{\boldsymbol{x}}_i).
\end{equation}
are derived from a quadrature rule on $C$, composed of nodes $\hat{\boldsymbol{x}}_i$ and weights $\hat{w}_i$.

Since the mappings $\mappingsB$ are bounded by the transformed level set $\tilde{\beta}$,
evaluating the height function $b^0$ of $B \in \mappingsB$,
\begin{equation}
    b^0 : \tilde{\beta}(\Pi_h(\tilde{\boldsymbol{x}}) + b^0(\tilde{\boldsymbol{x}}) \boldsymbol{e}_h) = 0 
\end{equation}
defined by the level set's zero-isocontour requires the Jacobian matrix $DA$.
The Jacobian matrix is required in the newton method
\begin{equation}
    \quad
    b^0_{n+1} = b^0_n - \frac{\tilde{\beta}(\tilde{\boldsymbol{x}}_n)}{\partial_{b^0} \tilde{\beta}(\tilde{\boldsymbol{x}}_n) }
    = b^0_n - \frac{\beta(A(\tilde{\boldsymbol{x}}_n))}{ \nabla \beta(A(\tilde{\boldsymbol{x}}_n)) \, D A \cdot \boldsymbol{e}_h } ,
    \quad
    \boldsymbol{x}_n = \Pi_h(\tilde{\boldsymbol{x}}) + b^0_n(\tilde{\boldsymbol{x}}) \boldsymbol{e}_h
    ,
\end{equation}
where the composition $\tilde{\beta}$ is differentiated.

\subsection{Surface Integrals}

Let the set $A_{\Gamma}\subseteq \mappingsA$ group all mappings with a codomain $V_i$
adjacent to the surface $\Gamma_{\alpha}$, so that $\Gamma_{\alpha} \cap V_i\neq \emptyset$. 
Moreover, let the mapping $E_A$ embed a unit square in the face 
of the domain of the mapping $A\in A_{\Gamma}$ that maps to $\Gamma_{\alpha}$.
More precisely, if we assume without loss of generality that 
$\Gamma_{\alpha}$ forms the upper face of the codomain of $A$ and $h_3 = z$, 
\begin{equation}
    A(C^{\uparrow}) \subseteq \Gamma_{\alpha}
    ,\quad
    C^{\uparrow} =\text{face}_z^{\uparrow}(C)
    ,
\end{equation}
then $E_A$ is the function, 
\begin{equation}
    E_A : C^{2} \rightarrow C^{\uparrow}, \quad E_A(x,y) = (x,y,1)^T,
\end{equation}
that embedds a unit square $C^2$ on the upper face of $C$.

Facilitated by $E_A$, 
the composition of mappings $T$,
\begin{equation}
    T: \mathbb{R}^2 \rightarrow \mathbb{R}^3, \quad T = A \circ E_A,
\end{equation}
maps a two-dimensional unit square onto $\Gamma_{\alpha}$. 
For each mapping $A \in A_{\Gamma}$,
the composition transforms the unit square to the surface integral,
\begin{equation}
    \oint_{\Gamma_{\alpha}} g(\boldsymbol{x}) \,\text{d}S  
    = \sum_{A \in A_{\Gamma}} \int_{C^2} g(T(\tilde{\boldsymbol{x}})) \, G_T(\tilde{\boldsymbol{x}})\,\text{d}S 
    \approx \sum_{A \in A_{\Gamma}} \sum_i w_{T,i} g(\boldsymbol{x}_{T,i}).
\end{equation}
Subsequently, the surface quadrature nodes $\boldsymbol{x}_{T,i}$ and weights $ w_{T,i}$ are derived, 
\begin{equation}
    w_{T,i} = G_T(\tilde{\boldsymbol{x}}_i) \tilde{w}_i,\quad  \boldsymbol{x}_{T,i} = T(\tilde{\boldsymbol{x}}_i),
\end{equation}
from a set of quadrature nodes $\tilde{\boldsymbol{x}}_i$ and weights $\tilde{w}_i$ on $C^2$.
Since $T$ is an embedding, the Gram determinant which generalizes the Jacobian determinant,
\begin{equation}
    G_T(\tilde{\boldsymbol{x}}) = \sqrt{\det\left(DT(\tilde{\boldsymbol{x}})^T \, DT(\tilde{\boldsymbol{x}}) \right)},
\end{equation}
encodes the geometry of the surface.

The integral on a surface of a domain confined by two level sets can be found 
by injecting the embedding $E_A$ into the composition of mappings $A \circ B$ of the volume integral.
We compose $T$, 
\begin{equation}
    T : \mathbb{R}^2 \rightarrow \mathbb{R}^3, T = A \circ E_A \circ B,
\end{equation}
and transform $\gamma_{\alpha}$ to a set of unit squares, 
\begin{equation}
    \oint_{\gamma_{\alpha}} g(\boldsymbol{x}) \,\text{d}S  
    = \sum_{A \in A_{\Gamma}} \sum_{B \in B_{\Omega}}\int_{C^2} g(T(\hat{\boldsymbol{x}})) \, G_{T}(\hat{\boldsymbol{x}})\,\text{d}S 
    \approx \sum_{A \in A_{\Gamma}} \sum_{B \in B_{\Omega}} \sum_i w_{T,i} g(\boldsymbol{x}_{T,i}),
\end{equation}
to construct a quadrature rule over $\gamma_{\alpha}$. 
The composition $T$ transforms the integral in two steps. 
First, the mapping $A \circ E_A$ maps a unit square onto a subset of the surface $\Gamma_{\alpha}$.
Then, after transforming $\tilde{\beta} = \beta(A \circ E_A)$ to this unit square, 
the transformed integral is treated as a two-dimensional volume integral. 
Consequently, the codomains of the mappings $B_{\Omega}$ tesselate the two-dimensional domain $\Omega_{\tilde{\beta}}$.

Similiarly, 
a quadrature rule over $\gamma_{\beta}$ 
is provided by the composition $T$,
\begin{equation}
    T : \mathbb{R}^2 \rightarrow \mathbb{R}^3, T = A \circ B \circ E_B.
\end{equation}
The composition restricts the domain of integration to $\Omega_{\alpha}$ and 
then maps the surface defined by $\tilde{\beta}$ to $C^2$,
\begin{equation}
    \oint_{\gamma_{\beta}} g(\boldsymbol{x}) \,\text{d}S  
    = \sum_{A \in \mappingsA} \sum_{B \in B_{\Gamma}} \int_{C^2} g(T(\tilde{\boldsymbol{x}})) \, G_{T}(\tilde{\boldsymbol{x}})\,\text{d}S 
    \approx \sum_{A \in \mappingsA} \sum_{B \in B_{\Gamma}} \sum_i w_{T,i} g(\boldsymbol{x}_{T,i}),
\end{equation}
using the mappings $E_B$ of each mapping $B \in B_{\Gamma}$. Corresponding to the definition of $A_\Gamma$, 
the set $B_{\Gamma}\subseteq B_{\Omega}$ 
groups the mappings with a codomain in $\Omega_{\tilde{\beta}}$ that is adjacent to $\tilde{\beta}$.

\subsection{Line Integrals}

Applying the composition of transformations $T$, 
\begin{equation}
    T: \mathbb{R} \rightarrow \mathbb{R}^3, T = A \circ E_A \circ B \circ E_B,
\end{equation}
yields the line quadrature rule, 
\begin{equation}
    \fint_{\partial \gamma} g(\boldsymbol{x}) \,\text{d}L  
    = \sum_{A \in A_{\Gamma}} \sum_{B \in B_{\Gamma}} \int_{C^1} g(T(\hat{\boldsymbol{x}})) \, G_{T}(\hat{\boldsymbol{x}})\,\text{d}L
    \approx \sum_{A \in A_{\Gamma}} \sum_{B \in B_{\Gamma}} \sum_i w_{T,i} g(\boldsymbol{x}_{T,i}).
\end{equation}
The composition transforms the integral in two steps. 
First, the mapping $A \circ E_A$ maps a unit square $C^2$ onto the surface $\Gamma_{\alpha}$.
Then, $B \circ E_B$ maps a unit line $C^1$ to the surface
defined by $\tilde{\beta} = \beta(A \circ E_A)$ in $C^2$.

The line quadrature nodes $\boldsymbol{x}_{T,i}$ and weights $ w_{T,i}$, 
\begin{equation}
    w_{T,i} = G_T(\tilde{\boldsymbol{x}}_i) \tilde{w}_i,\quad  \boldsymbol{x}_{T,i} = T(\tilde{\boldsymbol{x}}_i),
\end{equation}
are derived from a one-dimensional set of quadrature nodes $\tilde{\boldsymbol{x}}_i$ and weights $\tilde{w}_i$ on $C^1$.

\section{Experiments} \label{section:experiments}

We conduct five numerical experiments to investigate the method\cite{Data}.
In the first experiment, we evaluate the edge length, surface area and volume of a domain with a sharp, oscillating edge.
Offering a curved topology with a Gaussian curvature of zero, the domain's shape is defined by trigonometric level sets. 
In the second experiment we investigate a shape with non-zero Gaussian curvature
and evaluate the edge length, surface area and volume of a lens defined by two intersecting spheres.
Examining both an intricate shape and integrand in the third experiment, 
we evaluate the line and surface integral of a trigonometric function on a curved toric section.
In the last two experiments, 
we investigate shapes described by a single level set that require a great amount of quadrature nodes 
due to subdivisions. 
We reduce the amount of quadrature nodes by inserting a second level set.

In each experiment we measure the error $e$,
\begin{equation}
    e = \left| G - \sum_i^{c^3} \sum_{T \in T_i} \sum_j^{n^d} w_{T,j} g\left(\boldsymbol{x}_{T,j}\right) \right| 
    \quad G = \fint_{\partial \gamma} g(\boldsymbol{x}) \,\text{d}L ,
    \oint_{\gamma_{\beta}} g(\boldsymbol{x}) \,\text{d}S ,
    \int_{\Omega} g(\boldsymbol{x}) \,\text{d}V ,
\end{equation}
which is the difference between the exact integral and its approximation. 
For the sake of brevity, $e$ denotes the error for every domain of integration, 
regardless if its a volume ($d=3$) a surface ($d=2$) or a line ($d=1$).
Sharing a basic setting, each experiment is conducted in a cube $K = [-1,1]^3$  
which is subdivided into $c^3$ uniform cells. 
To compute the quadrature rules for each cell, 
we create a set of mappings $T_i$ for each $i$-th cell.
Subsequently each mapping $T$ of $T_i$ is applied to transform a tensorized Gauss quadrature rule constisting of $n^d$ quadrature nodes,
constructing the transformed quadrature nodes, $\boldsymbol{x}_{T,j}$, and weights, $w_{T,j}$.

\subsection{Oscillating Edge}
\label{subsection:Oscillating}

Let the level sets $\alpha$ and $\beta$,  
\begin{gather}
    \alpha(x,y,z) = z - \frac{1}{5} \sin\left(\frac{20 \pi x}{11}\right), \\
    \beta(x,y,z) = y - \frac{1}{5} \sin\left(\frac{20 \pi x}{11}\right), \\
\end{gather}
define a trigonometric domain in the cube $K$.
The intersection of the level sets creates a sharp edge, oscillating around the $x$-axis. 
Focusing on the geometry of the domain of integration, we integrate a constant function,
\begin{equation}
    g(x,y,z) = 1,
\end{equation}
and evaluate the length of $\partial \gamma$, the surface area of $\gamma_{\beta}$ and the volume of $\Omega$.
Choosing a constant $g$ reduces the integrand of each transformed integral to
the Jacobian determinant and removes any influence of $g$.

To compute the exact solution, 
we precisely evaluate the integrals with an established quadrature method~\cite{Shampine}, 
\begin{gather}
    \fint_{\partial \gamma} 1 \,\text{d}L = 2.9018098242473137628716230441128...\,, \\
    \oint_{\gamma_{\beta}} 1 \,\text{d}S = 2.5048230500093248969863804012397...\,, \\
    \int_{\Omega} 1 \,\text{d}V  = 2.0431849934260147426243415665995... \, ,
\end{gather}
using the analytic expressions of $\gamma_{\alpha}$ and $\gamma_{\beta}$.
 
\begin{figure}
    \centering
    \includegraphics{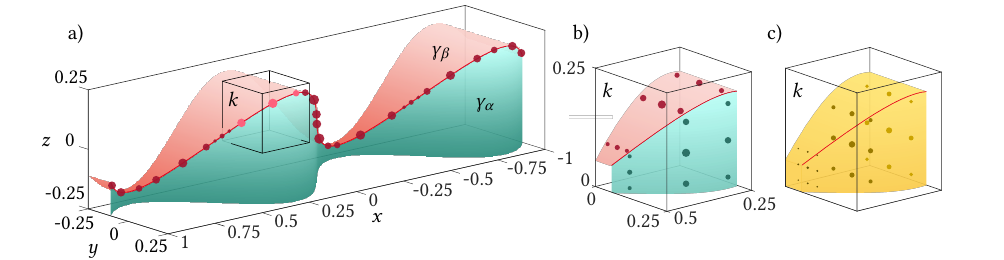}
    \caption{Two trigonometric surfaces $\gamma_{\alpha}$ and $\gamma_{\beta}$ that form an oscillating edge at their intersection. 
    Each filled circle denotes a node of a transformed Gauss quadrature rule, $n=3$. 
    The radius of each node is proportional to its weight.
    The left picture shows the quadrature rule obtained for the oscillating edge 
    on a grid consisting of $c^3=8^3$ uniform cells (a).
    The transformed quadrature rules on the surfaces (b) and in the volume (c) are depicted for a single cell $k$ on the right.}
    \label{figure:OscillatingLine}
\end{figure}

Figure~\ref{figure:OscillatingLine}
depicts a section of the domains of integration $\partial \gamma, \gamma_{\alpha},\gamma_{\beta}$ and $\Omega$ in proximity of $\partial \gamma$
and the corresponding quadrature rules for $n = 3$ and $c = 8$. 
The nodes of the line quadrature rule of $\partial \gamma$ are located on the edge, 
showing the characteristic spacing of the underlying Gauss rule.
Groups of $n$ nodes are positioned with varying distance, 
depending on the length of the intersection of the edge and the cell. 
By construction, each intersection contains at least $n$ nodes.
Similarly, the nodes of the surface quadrature rules of $\gamma_{\alpha},\gamma_{\beta}$ are located on the surface 
and the nodes of the volume quadrature rule  
of $\Omega$ are located in the volume.

\begin{figure}
    \centering
    \includegraphics{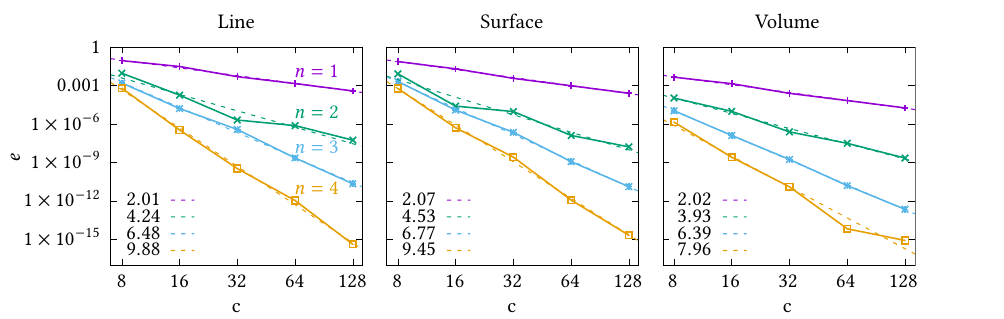}
    \caption{The error $e$ obtained during the evaluation of the integrals of the oscillating edge.
    Each curve is associated to $n$ underlying Gauss nodes. 
    When the number of cells $c$ per coordinate axis is increased, 
    the error converges with the order depicted in the bottom left of each graph.}
    \label{figure:OscillatingConvergence}
\end{figure}

If $c$ is successively increased, the observed convergence order 
is in good accordance with the expected convergence order of $2n$, Figure~\ref{figure:OscillatingConvergence} shows.
It shows the error $e$ and the order of convergence obtained while 
evaluating the line, surface and volume integrals. 

\subsection{Spherical Lens}

\label{subsection:Intersecting}
Creating a quarter lens with a sharp edge,
let the zero-isocontours of $\alpha$ and $\beta$ describe two intersecting spheres,  
\begin{gather}
    \alpha(x,y,z) = (x+1)^2 + (y+1)^2 + \left(z + \frac{49}{100}\right)^2 - \left(\frac{9}{10}\right)^2, \\
    \beta(x,y,z) = (x+1)^2 + (y+1)^2 + \left(z - \frac{51}{100}\right)^2 - \left(\frac{9}{10}\right)^2, \\
    K = [-1,1]^3,
\end{gather}
in cube $K$.

Focusing on the geometry of the domain of integration again, 
we evaluate the length of $\partial \gamma$, 
the surface area of $\gamma_{\beta}$ and the volume of $\Omega$ 
and choose $g$ accordingly,  
\begin{equation}
    g(x,y,z) = 1.
\end{equation}
The exact edge length, surface area and volume,  
\begin{gather}
    \fint_{\partial \gamma} 1 \,\text{d}L = \frac{\sqrt{14}}{10}\pi, \\
    \oint_{\gamma_{\beta}} 1 \,\text{d}S = \frac{9}{50}\pi, \\
    \int_{\Omega} 1 \,\text{d}V  = \frac{23}{375}\pi,
\end{gather}
can be found analytically.

\begin{figure}
    \centering
    \includegraphics{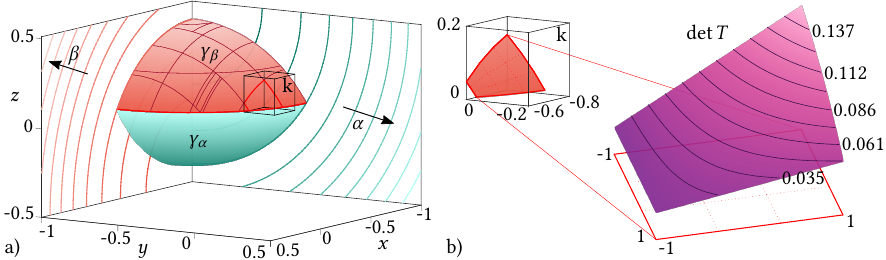}
    \caption{A lens created by intersecting spheres, defined by level sets $\alpha$ and $\beta$ (a). 
    The red lines on the surface $\gamma_{\beta}$ outline the codomains of each surface mapping 
    created on a grid consisting of $c^3=10^3$ uniform cells. 
    The right side (b) shows a detailed view of the Jacobian determinant, $\det T$, of a surface mapping
    belonging to a surface patch bounded by $\partial \gamma$ located in cell $k$. 
    }
    \label{figure:ufo}
\end{figure}

Showing the surfaces $\gamma_{\alpha}$ and  $\gamma_{\beta}$ of the lens, 
Figure~\ref{figure:ufo} illustrates two properties of the surface mappings constructed on a grid with $c^3=10^3$ cells. 
First, it outlines the codomains of the surface mappings of $\gamma_{\alpha}$ of each 
cell with red lines.
The codomains tesselate $\gamma_{\alpha}$ and vary in shape, depending on the relative position of the cells and the surface.  
Second, the figure illustrates how the geometry of a surface patch in a cell containing a part of the edge
is encoded in the Jacobian determinant. 
Encoding the transformation to the unit square, the Jacobian determinant is non-polynomial.

\begin{figure}
    \centering
    \includegraphics{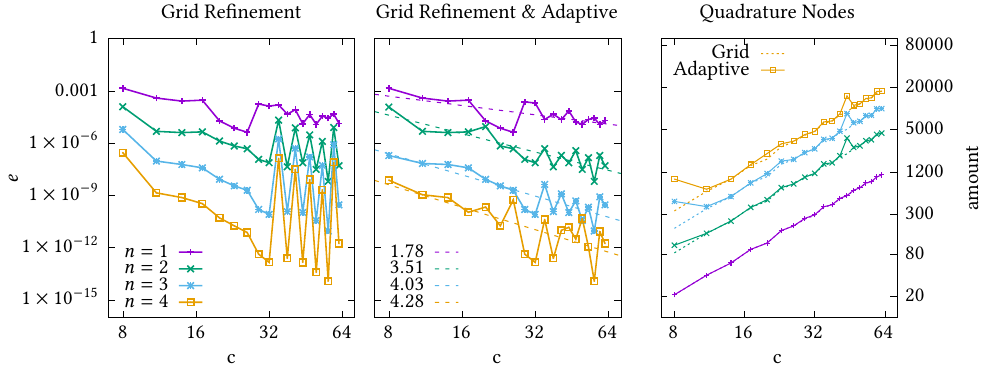}
    \caption{Surface area error of the spherical lens under grid refinement and grid refinement in combination with adaptive quadrature.
        A linear trend is indicated by a dotted line for each transformed Gauss rule $n$.
        The right panel shows the total amount of quadrature nodes of each rule. 
        The dotted lines depict the amount of nodes created while refining the grid
        and
        the solid lines depict the amount of nodes created by additionally applying adaptive quadrature.
        }
    \label{figure:UfoSurfaceConvergence}
\end{figure}

When the grid is refined, 
the convergence order of the surface error is as expected for $c$ below approximately 32, the left panel of
Figure~\ref{figure:UfoSurfaceConvergence} shows. 
After $c$ has reached approximately 32 however, the convergence order deteriorates 
and for some values of $c$ the error is significantly higher than the trend line.
This effect becomes more pronounced when $n$ is increased.
The fact that the error decreases when $n$ is increased and the value of $c$ is held constant,
suggests that
the irregular error is caused by the Jacobian determinant $\det T$ in the integrand.
Depending on its surrounding cell, 
the Jacobian determinant $\det T$ of a mapping may be a ill-behaved non-polynomial function.
When this integrand is badly approximated by the polynomial Gauss quadrature rule,
the error in these cells converges with increasing $n$, but lies over the average error level nonetheless.

To recover accuracy, we apply adaptive quadrature in the domains of the mappings.
Since the domain of each mapping is a $i$-dimensional hypercube $C$, 
many methods can be applied, such as Gauss-Kronrod quadrature or methods specially tailored to rational functions.
For simplicity, we apply a straightforward adaptive quadrature method.
The quadrature rule $Q$ belonging to the mapping is constructed by placing tensorized Gauss nodes and weights in $C$, 
so that $|Q| = n^i$.
If the estimated error $e_Q$ of the quadrature rule is greater than a threshold $\tau$,
\begin{equation}
    e_Q > \tau,
\end{equation}
$C$ is subdivided into $2^i$ hypercubes of equal size.
Repeating the process, the quadrature rule $\tilde{Q}$ belonging to the mapping is constructed 
by placing tensorized Gauss nodes and weights in every new hypercube, so that $|\tilde{Q}| = 2^i n^i$.
If the estimated error in one of the new hypercubes is too great still, it is subdivided again. 
And so on until a predetermined depth is reached.

Error estimation is a critical component of an adaptive quadrature method 
and many error erstimation schemes have been presented~\cite{Gonnet}.
We choose a local relative error,
\begin{equation}
    e_Q = \left| 
        \frac{\sum_{i}  g\left( \boldsymbol{x}_{Q,i} \right) w_{Q,i} 
            - \sum_{i}  g\left(\tilde{\boldsymbol{x}}_{\tilde{Q},i}\right) \tilde{w}_{\tilde{Q},i}}
        {\sum_{i}  g\left(\boldsymbol{x}_{Q,i}\right) w_{Q,i}}
    \right|
\end{equation}
relating the result of quadrature rule $Q$ to the result of the subdivided quadrature rule $\tilde{Q}$. 

\begin{figure}
    \centering
    \includegraphics{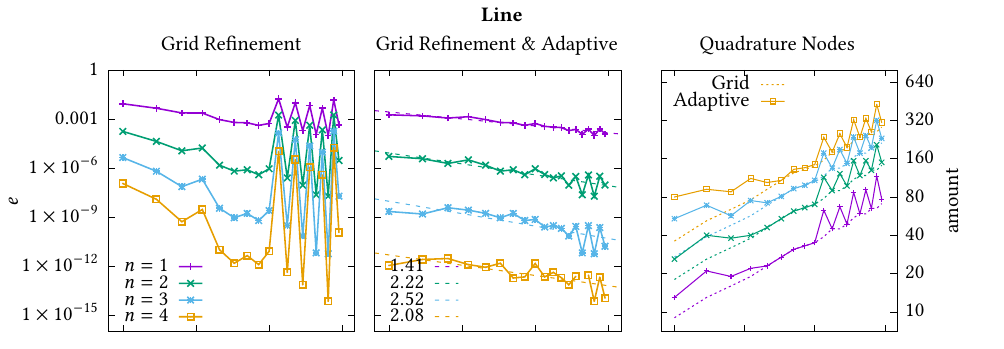}
    \includegraphics{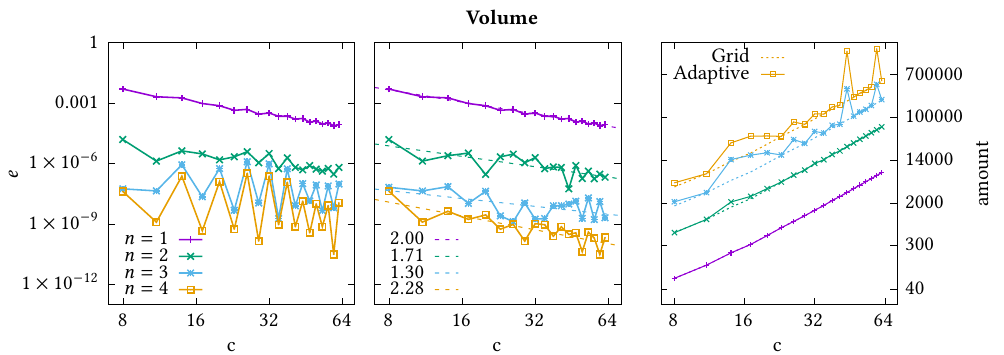}
    \caption{
        Line length error and volume error of the spherical lens under grid refinement and grid refinement 
        in combination with adaptive quadrature.
        A linear trend is indicated by a dotted line for each transformed Gauss rule $n$. 
        The right panel shows the total amount of quadrature nodes of each rule. 
        The dotted lines depict the amount of nodes created while refining the grid 
        and
        the solid lines depict the the amount of nodes created by additionally applying adaptive quadrature.
        }
    \label{figure:UfoVolumeLineConvergence}
\end{figure}

The middle panel of Figure~\ref{figure:UfoSurfaceConvergence} 
depicts the effect of adaptive quadrature on the surface integrals of the lens with the thresholds $\tau_n$ of each $n$. 
It lists the error for the values $\tau_1 = 0.1, \tau_2 = 0.001, \tau_3 = 1 \times 10^{-5}, \tau_4 = 1 \times 10^{-7}$,
showing that accuracy can be recovered with the adapative quadrature method. 
However, it decreases the error of lower $c$ to greater extent, 
flattening its trend and leading to a convergence order which is lower than $2n$. 
This is also reflected by 
the total number of quadrature nodes created
which is shown for each $n$ in the right panel of Figure~\ref{figure:UfoSurfaceConvergence}.
For lower values of $c$, the number of quadrature nodes of the adaptive quadrature method lies above the 
number of quadrature nodes of the method without adaptation. For higher values of $c$, 
the number of quadrature nodes created adaptively
slightly lies above the number of nodes created without adaptation in
the cases where the error of the method without adaptation produces a disproportionate error. 
Otherwise they are almost equal.

Deviation from the expected error at certain resolutions $c$ can also be observed in the computation of the edge length and the volume, 
illustrated by Figure~\ref{figure:UfoVolumeLineConvergence}.  
In each case, adaptive quadrature is applied as a remedy, effectively removing unwanted deviation of the error.
The thresholds $\tau_n$ of each $n$ used for adaptive quadrature of the line are,    
$\tau_1 = 0.002, \tau_2=1 \times 10^{-5}, \tau_3 = 1 \times 10^{-8}, \tau_4 = 1 \times 10^{-11} $.
The respective thresholds used for the volume are, $\tau_1 = 1, \tau_2 = 0.05, \tau_3 = 0.001, \tau_4 = 1 \times 10^{-4}$.
When choosing the threshold $\tau$ carefully, the amount of quadrature nodes increases only slightly, 
the right panel of Figure~\ref{figure:UfoSurfaceConvergence} and Figure~\ref{figure:UfoVolumeLineConvergence} shows.

\begin{figure}
    \centering
    \includegraphics{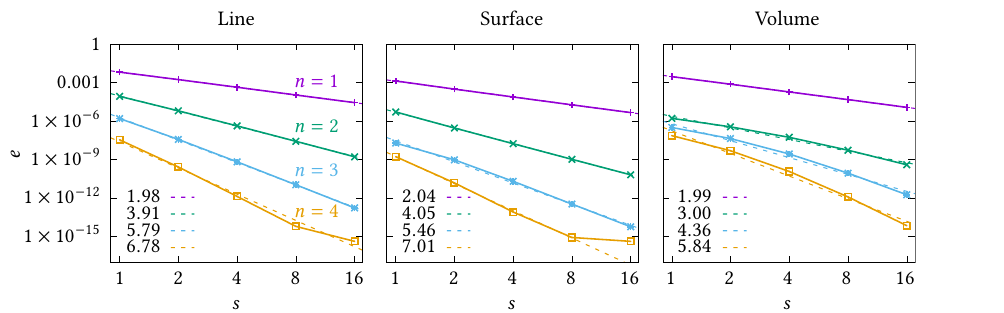}
    \caption{The error $e$ obtained during the evaluation of the integrals of the spherical lens 
    for different numbers of subdivisions~$s$ on a grid of $c = 10$ cells per coordinate axis.
    To increase accuracy, $s^i$ quadrature rules $n$ are placed in each $i$-dimensional domain of each mapping.
    The number in the bottom left of each graph denotes the order of convergence.
    }
    \label{Figure:SubdivisionError}
\end{figure}

The previous experiments show that the error decreases when increasing $n$ and holding $c$ constant.
Motivated by this, we choose a grid, $c=10$, 
and directly refine the quadrature rules in each domain $C$ of the grid's mappings. 
Each $i$-dimensional domain $C$ is subdivided evenly into $s^i$ hypercubes 
and a tensorized Gauss rule $n$ placed in each.
Providing the trend of the error when successively increasing $s$, 
Figure~\ref{Figure:SubdivisionError} suggests
an order of convergence of about $2n$ for each curve $n$.

\subsection{Toric Section}

\begin{figure}
    \centering
    \includegraphics{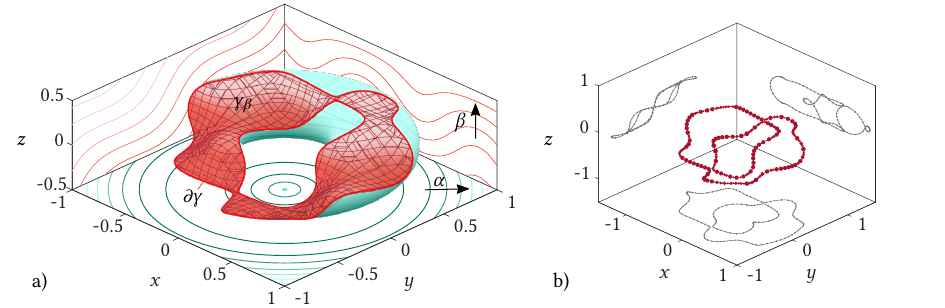}
    \caption{
        A tilted wavy surface intersects a torus, generating a curved toric section (a).
        The red lines on the surface $\gamma_{\beta}$ outline the
        codomains of each surface mapping created on a grid consisting of $c^3 = 30^3$ uniform cells.
        The right (b) shows a line quadrature rule on the the boundary line $\partial \gamma$ with $n = 1$. 
        The projections of the line onto the $xy$-, $xz$- and $yz$-plane are depicted in grey. 
        }
    \label{figure:Torus}
\end{figure}

Respectively describing a torus and a tilted, wavy surface, the level sets,  
\begin{gather}
    \alpha(x,y,z) = \left(\sqrt{x^2 +y^2} - \frac{3}{5} \right)^2 + z^2 -\left(\frac{3}{10}\right)^2, \\
    \beta(x,y,z) = -\frac{1}{2} y + \frac{7}{8} x - \frac{1}{5} \sin(\frac{20\pi}{11} x)  \sin\left(\frac{35\pi}{22}  y + \frac{10\pi}{11}  z \right) ,
\end{gather}
generate a curved surface $\gamma_{\beta}$ that is skewed to the grid (see Figure~\ref{figure:Torus}). 
The function $g$ is a combination of two trigonometric functions,
\begin{equation}
    g(x,y,z) = \frac{1}{5} \sin(\frac{20\pi}{11} x)  \sin\left(\frac{35\pi}{22}  y + \frac{10\pi}{11}  z \right).
\end{equation}
Since $\gamma_{\beta}$ and $g$ are line-symmetric with respect to the $x$-axis,
the integrals of $g$ along $\partial \gamma$ and over $\gamma_{\beta}$,
\begin{gather}
    \fint_{\partial \gamma} g(\boldsymbol{x}) \,\text{d}L = 0, \\
    \oint_{\gamma_{\beta}} g(\boldsymbol{x}) \,\text{d}S = 0,
\end{gather}
vanish.

The curved intersection $\partial \gamma$ of the torus and the wavy sheet is depicted in the right picture of Figure~\ref{figure:Torus}.  
Its intricate shape is outlined by
the quadrature nodes of a line quadrature rule. The line quadrature rule is
immersed in a grid of $c^3 = 30^3$ cells and consists of $n = 1$ Gauss nodes per nested mapping.
Although the shape of $\partial \gamma$ is complicated, the nodes are placed evenly along the boundary.
From the projections of $\partial \gamma$ onto the $xy$-, $xz$- and $yz$-plane Figure~\ref{figure:Torus} it can be seen, 
that the boundary is line-symmetric. 

\begin{figure}
    \centering
    \includegraphics{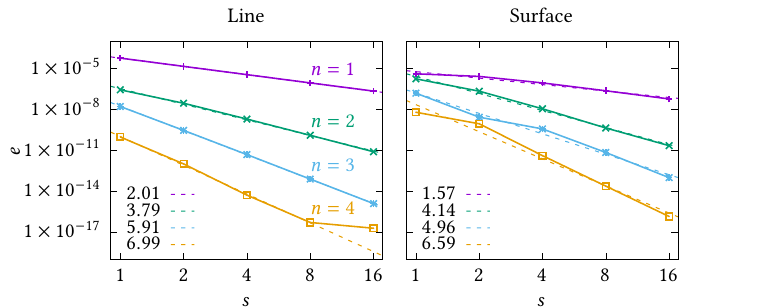}
    \caption{The error $e$ obtained during the evaluation of the line and surface integrals of the toric section
        with $n$ underlying Gauss nodes.
        The set of mappings created on a grid of $c^3 = 30^3$ cells is kept constant.
        To increase accuracy, $s^i$ quadrature rules $n$ are placed in each $i$-dimensional domain of each mapping.
        The number in the bottom left of each graph denotes  
        the order of convergence.
    }
    \label{figure:TorusConvergence}
\end{figure}

The function $g$ is a composition of trigonometric functions that creates a wavy, smooth field throughout $K$.
We evaluate the line and surface integral of $g$ with a quadrature rule with $n$ nodes on a grid of $c^3 = 30^3$ cells. 
When subdividing the quadrature rules direclty in each domain of the grid's mappings,
we measure a convergence rate of the error of about $2n$, as illustrated by Figure~\ref{figure:TorusConvergence}.
This shows that high order convergence can be achieved for intricate geometries and functions.

\subsection{Thin Cylindrical Sheet}

\label{subsection:cylindrical}

Depending on the relative position of the zero-isocontour and $K$, 
the construction of an exact nested mapping may require a high number of subdivisions, 
even for a quadratic level set.
For example, let the surface of a quadratic level set $\alpha$, 
\begin{gather}
    \alpha (x,y,z) = \tilde{\alpha}(x,y,z)^2 - l^2 \quad
    \tilde{\alpha}(x,y,z) = \left(\frac{4}{5}x - y\right) 
\end{gather}
be two parallel planar surfaces separated by distance $l$, creating a slanted sheet. 
Both components of the linear gradient of $\alpha$,
\begin{equation}
    \nabla \alpha =  2 \tilde{\alpha} \nabla \tilde{\alpha} = 2 \tilde{\alpha}
    \begin{pmatrix}
        \frac{4}{5}\\ 
        -1\\
        0
        \end{pmatrix}
\end{equation}
contain the zero-isocontour $\Gamma_{\tilde{\alpha}}$,
\begin{equation}
    \Gamma_{\tilde{\alpha}} = \left\{ \boldsymbol{x} \in K \,|\, \tilde{\alpha}(\boldsymbol{x} ) = 0  \right\}
\end{equation} 
which is pinched between the two planar surfaces in $\Gamma_{\alpha}$.

\begin{figure}
    \centering
    \begin{minipage}{0.4\textwidth}
        \includegraphics{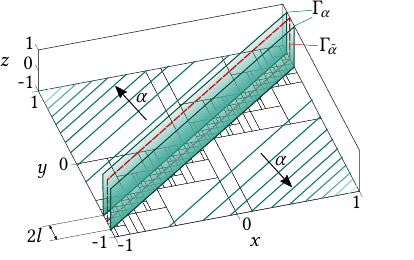}
    \end{minipage}
    \begin{minipage}{0.4\textwidth}
        \begin{tabular}{r l l}
            \toprule
            $l$ & Subdivision & Contour Splitting\\
            \midrule
            0.1 & 2426 & 2\\
            0.03 & 29010 & 2 \\
            0.01 & 248162 & 2 \\
            \bottomrule
            \\ [-1em]
            \multicolumn{3}{c}{ Amount of volume quadrature nodes }  \\
        \end{tabular}
    \end{minipage}
    \caption{Slanted sheet created by two parallel planar surfaces separated by distance $l$ of a single level set $\alpha$.
        Pinched between the surfaces, the zero-isocontour of both gradients $\Gamma_{\tilde{\alpha}}$ is indicated by a dashed red line.
        The subsets of $K$ created by Algorithm~\ref{alg:graph} with unrestricted subdivision are outlined in the $xy$-Plane.
        Comparing the subdivision and contour splitting approach, 
        the table on the right lists the number of nodes contained in volume quadrature rules, $n = 1$, for different values of $l$.
   }
    \label{figure:slantedSheet}
\end{figure}

When Algorithm~\ref{alg:graph} tries to represent the surfaces by a graph,
it will select a height direction $h$
and subdivide $K$ until the zero-isocontour of $\partial_h \alpha$  
exclusively lies in subsets of $K$ not containing $\gamma_{\alpha}$ (see Figure~\ref{figure:slantedSheet}).
Since the the zero-isocontour of $\partial_h \alpha$, $\Gamma_{\tilde{\alpha}}$, is pinched between the two surfaces of $\Gamma_{\alpha}$, 
decreasing distance $l$ causes a substantial increase of subdivisions during the construction of quadrature rules.
If however $\Omega_\alpha$ is split along $\Gamma_{\tilde{\alpha}}$ before creating the quadrature rule, 
subdivisions are not required and the number of nodes constructed is reduced greatly.
This is demonstrated by the table on the right side of Figure~\ref{figure:slantedSheet} which shows
the number of nodes created for a volume integral for different values of $l$.
In the case of $l=0.01$, splitting along the zero-isocontour leads to approximately $62^3$ fewer nodes. 

Reducing the amount of quadrature nodes is realized in two steps. 
First, we split $\Omega_{\alpha}$ along the zero-isocontour of $\partial_h \alpha$ into $\Omega'_{\alpha}$ and $\Omega''_{\alpha}$,
\begin{equation}
    \Omega_{\alpha} = \Omega'_{\alpha} \dot \cup \Omega''_{\alpha},
\end{equation}
after the algorithm has chosen a height direction $h$. 
The two new domains result from
intersecting $\Omega_{\alpha}$ with the domains $\Omega_{\beta}$ and $\Omega_{\beta^-}$,
\begin{equation}
    \Omega'_{\alpha} = \Omega_{\alpha} \cap \Omega_{\beta}, 
    \quad \Omega''_{\alpha} = \Omega_{\alpha} \cap \Omega_{\beta^-},
    \quad \beta = \partial_h \alpha,
    \quad \beta^- = -\partial_h \alpha.
\end{equation}
defined by the level sets $\beta$ and $\beta^-$.
Second, we construct quadrature rules on the two new domains, applying the approach for two level sets.
Because the zero-isocontour of $\partial_h \alpha$ lies on their boundary, 
subdivisions are not required on the new domains. 

\begin{figure}
    \centering
    \begin{minipage}{0.575\textwidth}
        \includegraphics{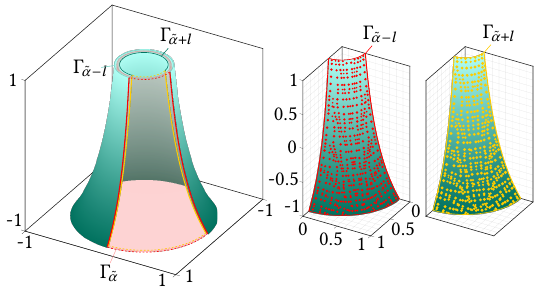}
    \end{minipage}
    \begin{minipage}{0.4\textwidth}
        \includegraphics{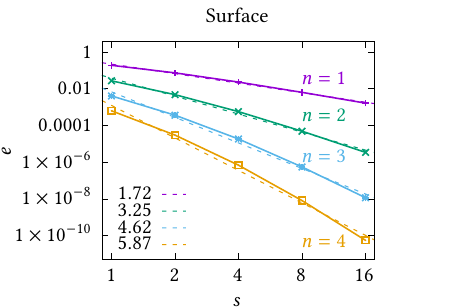}
    \end{minipage}
    \caption{A thin cylindrical sheet with two concentric surfaces $\Gamma_{\tilde{\alpha}-l}$ and $\Gamma_{\tilde{\alpha}+l}$, 
        encasing a zero-isocontour $\Gamma_{\tilde{\alpha}}$ of the level set's gradient drawn as a red, transparent surface.  
        A section of the surface quadrature nodes, $n=1$, of both surfaces is displayed in the middle.
        The quadrature nodes are created on a grid containing $c^3 = 20^3$ cells, indicated by grey lines.
        A convergence plot of the error of the surface area is shown on the right, 
        where $s^i$ quadrature rules $n$
        are placed in each $i$-dimensional domain of each mapping.
        }
    \label{figure:volcano}
\end{figure}

This ansatz can be applied generally to domains of smooth levels sets 
with derivatives unfavorable for subdivision. 
For example, let $\alpha$ be a thin cylindrical sheet, 
\begin{equation}
    \alpha (x,y,z) = \tilde{\alpha}(x,y,z)^2 - l^2 \quad
    \tilde{\alpha}(x,y,z) = x^2 + y^2 + \left(z + \frac{6}{5}\right)^{-2} \quad 
    l = \frac{1}{50}
\end{equation}
that has two concentric surfaces $\Gamma_{\tilde{\alpha}-l}$ and $\Gamma_{\tilde{\alpha}+ l}$,
\begin{equation}
    \Gamma_{\tilde{\alpha}-l} = \{ \boldsymbol{x} \in K \,|\, \tilde{\alpha}(\boldsymbol{x} )-l  = 0  \}, \quad
    \Gamma_{\tilde{\alpha}+ l}= \{ \boldsymbol{x} \in K \,|\, \tilde{\alpha}(\boldsymbol{x} )+l  = 0  \},
\end{equation}
that are curved and cylindrical. 
The gradient of $\alpha$ contains the zero-isocontour $\Gamma_{\tilde{\alpha}}$ of the curved cylinder $\tilde{\alpha}$,
\begin{equation}
    \nabla \alpha(x,y,z) = 2\tilde{\alpha} \, \nabla \tilde{\alpha},
\end{equation}
which lies between $\Gamma_{\tilde{\alpha}-l}$ and $\Gamma_{\tilde{\alpha}+l} $ (see Figure~\ref{figure:volcano}).

By first splitting the domain along $\Gamma_{\tilde{\alpha}}$, 
quadrature rules for the thin sheet defined by a single level set can be found without excessive subdivision.
The middle panel of Figure~\ref{figure:volcano} shows the quadrature nodes of both surfaces of a quadrature rule with $n=1$. 
They are obtained in a grid consisting of $c^3 = 20^3$ cells.
The amount of quadrature nodes depicted in the figure is approximately proportional to the number of cells intersecting the surfaces.  
This reflects a low number of subdivisions, 
although both surfaces are close and intersect the same cells.

In order to measure the order of convergence, 
we evaluate the surface areas of $\Gamma_{\tilde{\alpha}-l}$ and $\Gamma_{\tilde{\alpha}+ l}$, 
\begin{equation}
    g(x,y,z) = 1,
\end{equation}
by applying the splitting ansatz. 
The exact surface areas can be found analytically, 
\begin{gather}
    \oint_{\Gamma_{\tilde{\alpha}} - l} 1 \, \text{d}S = 6.775163182554902237379363684639... ,\\ 
    \oint_{\Gamma_{\tilde{\alpha}}+ l} 1 \, \text{d}S = 6.2377886792891965132267781181652... ,
\end{gather}
using that the surfaces are a revolution arround the $z$-axis.
We create the nested mappings required to create the quadrature rule in a grid consisting of $c^3 = 20^3$ cells 
and successively subdivide the domains of each nested mapping.
The resulting error can be seen in the right panel of Figure~\ref{figure:volcano}.
Despite the complex geometry of the integral's domain,  
the error obtained for $n=1,2,3,4$ converges exponentially with an order of approximately $1.5n$.

\subsection{Wavy Cylinder}

If evaluation of $g$ is computationally expensive, constructing a minimal amount of quadrature nodes is advantageous. 
Unfortunately, not only thin sheets immersed in a coarse grid can cause 
a high number of subdivisions leading to many quadrature nodes,
like it is the case in subsection~\ref{subsection:cylindrical}.
Even well resolved bodies can cause a high number of subdivisions when exact mappings are required.

\begin{figure}
    \centering
    \includegraphics{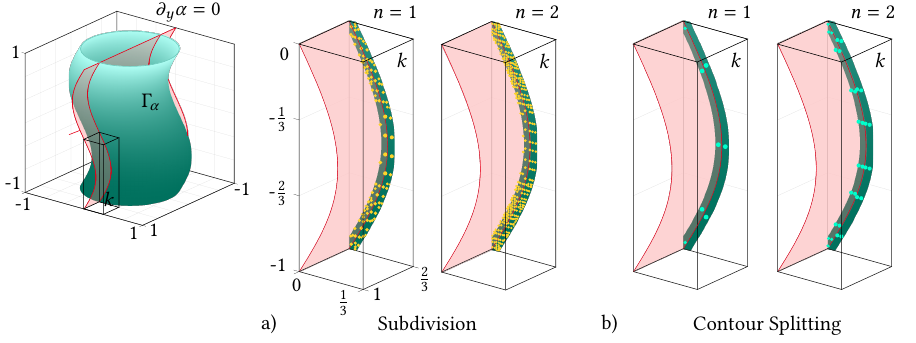}
    \caption{A wavy cylinder in a grid consisting of  $c^3 = 6^3$ cells that is slightly shifted.
        Each cell of the grid is indicated by grey lines.
        The two panels display the surface quadrature nodes of a subset of cells $k$ of the grid. 
        The surface quadrature nodes are constructed by (a) using subdivision and (b) splitting the cells along the zero-isocontour of $\partial_y \alpha$
        which is indicated by a red, transparent surface.
        }
    \label{figure:wavyCylinder}
\end{figure}

Let $\alpha$ describe a wavy cylinder, 
\begin{equation}
    \alpha(x,y,z) = x^2 + \left(y + \frac{1}{5} \sin(\pi z) \right)^2 - \left(\frac{2}{3} + l\right)^2,
\end{equation}
that is immersed in a grid of $c^3 = 6^3$ cells. 
A small perturbation $l$ of the radius leads to a high number of quadrature nodes
due to subdivision
in the cells located at the intersection of the zero-isocontour $\Gamma_{\tilde{\alpha}}$ of $\partial_y \alpha$, 
\begin{equation}
    \partial_y \alpha (x,y,z) = 2 \left(y + \frac{1}{5} \sin(\pi z) \right),
\end{equation}
and $\Gamma_\alpha$, Figure~\ref{figure:wavyCylinder} shows for $l= 0.001$. 

As described in subsection~\ref{subsection:cylindrical}, this is caused when $\Gamma_{\tilde{\alpha}}$
is pinched between two subsets of $\Gamma_{\alpha}$. 
For example, let us look at the three cells $k$ that 
are located at the intersection of $\Gamma_{\tilde{\alpha}}$ and $\Gamma_\alpha$ depicted in Figure~\ref{figure:wavyCylinder}.
Due to the slight shift of the wavy cylinder in $x$-direction, $k$ contains a thin region of $\Omega_{\alpha}$.
On the intersection of $\Omega_{\alpha}$ and the face $p^\downarrow = \text{face}^{\downarrow}_x(p)$ of the cell $p \in k$, 
the line $\Gamma_{\tilde{\alpha}} \cap p^\downarrow$ is pinched between two close surface sections $\Gamma_{\alpha} \cap p^\downarrow$.
When Algorithm~\ref{alg:graph} constructs the graph on $p$, 
it successively subdivides $p$ until the face $p'^{\downarrow}$ of each subcell $p'$
intersects only one of the three lines, 
\begin{equation}
    p'^{\downarrow} = \text{face}^{\downarrow}_x(p'), \quad
    p'^{\downarrow} \cap \Gamma_{\tilde{\alpha}} = \emptyset \,  \lor \, p'^{\downarrow} \cap \Gamma_{\alpha} = \emptyset,
\end{equation}
resulting in a high number of subdivisions.

To decrease the amount of subdivisions required to find exact mappings, 
we split the cells along the zero-isocontour of $\partial_y \alpha$ like in subsection~\ref{subsection:cylindrical}.
Figure~\ref{figure:wavyCylinder} shows the surface quadrature nodes 
in the three cells $k$ for both the approach with subdividing and the approach with contour splitting. 
The two sets of surface quadrature nodes displayed in the figure are constructed for the same distance $l = 0.001$ with $n=1$ and $n=2$. 
By comparing the surface quadrature nodes constructed by both approaches, 
it can be seen that contour splitting results in substantially fewer surface quadrature nodes for both values of $n$.
It is clear from the results of subsection~\ref{subsection:cylindrical}, 
that the difference in the amount of quadrature nodes of the two approaches further increases, 
when $l$ is decreased.

\section{Summary and Outlook} \label{section:summary}
In this paper we present a high-order method
that provides numerical integration of smooth functions on volumes, surfaces, and lines defined implicitly by two intersecting level sets.
Our approach allows for general smooth level sets 
and can map arbitrary quadrature rules defined on hypercubes to the curved domains of the integrals. 
Conceptually, the mappings required to create the quadrature rules 
are constructed by composing two nested mappings. 
The first nested mapping maps a hypercube to the domain defined by the first level set.
The second nested mapping maps a hypercube to the domain defined by
the second level set function composed with the first nested mapping which mirrors the domain of the integral.
Algorithmically, each nested mapping is generated by recasting the zero-isocontours of the level sets as the graphs of height functions.
This is achieved by augmenting a method developed for the integration on domains bounded by single level sets\cite{Saye}.

The practicality and usefullness of this approach is demonstrated in various numerical experiments.
We apply the approach to compute the volume, the surface area, and the intersection line length
of an oscillating edge formed by two trigonometric level sets and
show high-order convergence of the error.
In the investigation of the same properties of a spherical lens defined by two quadratic level sets, 
we measure high-order convergence 
when the quadrature rules are refined on a constant grid.
Fluctuating error levels introduced by the relative position of the grid and the domain of integration are successfully 
removed by applying an adaptive quadrature scheme, recovering high accuracy and 
showing that adapative quadrature conveniently suits the method.
Finally, to demonstrate high-order convergence of the error of an integral
on a complex non-polynomial geometry with an intricate integrand, 
a trigonometric function is integrated on a toric section.

In addition, two experiments that target the internal subdivision scheme of the method are conducted.
In the first experiment, a level set is inserted between the surfaces of a thin cylindrical sheet 
to construct a high-order surface quadrature rule without requiring internal subdivisions. 
The thin cylindrical sheet is defined by a rational level set and is
immersed in a coarse grid which otherwise would require the method to subdivide excessively.
In practical applications, 
internal subdivisions should be avoided because they
cause a high number of quadrature nodes which increases computational cost. 
Moreover, they reduce accuracy
when a low-order fallback method is applied to limit the subdivision depth.
In the second experiment,
a level set is added to reduce 
the internal subdivisions required for integration on a wavy cylinder. 
This shows that integrating on domains defined by two level sets can be advantageous, even for simple geometries.  

Beyond the close investigation of the method's general computational cost,
there are two future directions important for practical applications.
First, since the method draws from quadrature methods defined 
on hypercubes, a wide range of methods for adaptive subdivision could be applied 
to create quadrature rules taylored to specific integrands.  
Adaptive quadrature enables quadrature rules with a preset accuracy
which often is more pratical than providing a certain convergence order. 
Second, the dependency on low-order fallback methods could be removed 
by inserting level sets algorithmically. 
Especially when integration is performed many times on a static domain, the reduced number of
quadrature nodes could lead to a reduction of computation time.   

Another promising goal is to construct mappings which add minimal difficulty to the transformed integral. 
So far, each mapping is created by linearly interpolating between the boundaries of the integral's domain.
When the quadrature rule is subsequently mapped to the domain of the integral,
the additional term appearing in the integral can be unnecessarily complex.
This could be mitigated by choosing a different interpolation approach 
to increase accuracy while keeping the same number of quadrature nodes.

In an upcoming publication, the method of this paper will be combined with an extended discontinuous Galerkin method to 
investigate multi-phase problems with contact lines by the authors. 
Accurately modeling the complex behaviour of dynamic contact lines is a problem of ongoing research
which could benefit from precisely resolved geometries.

\subsection{Acknowledgements}
This research was funded by Deutsche Forschungsgemeinschaft (DFG, German Research Foundation) - 422800359,
and supported by the Graduate School CE within the Centre for Computational Engineering at TU Darmstadt.

\newpage 
\appendix
\section{Hessian} 
The hessian $H_{\tilde{\beta}}(\tilde{x}, \tilde{y},\tilde{z})$  of the composition $\tilde{\beta}(\tilde{\boldsymbol{x}}) = \beta(A(\tilde{\boldsymbol{x}}))$,
\begin{equation}
    H_{\tilde{\beta}}(\tilde{\boldsymbol{x}}) = H_{\beta}(A(\tilde{\boldsymbol{x}})) = D A^T \, H_{\beta} \, D A
    + \sum_{i \in \{x,y,z\}} \partial_{i} \beta H_{i},
\end{equation}
requires the hessians of $H_i$ of $A$.

Exhibiting increasing complexity, the Hessians $H_{x},H_{y}$ and $H_{z}$ reflect the recursive structure of the nested mapping $M$.
While the Hessian $H_{x}$ is composed of zeros, 
\begin{equation}
    H_{x} = 
    \begin{pmatrix}
        0 & 0 & 0 \\
        0 & 0 & 0 \\
        0 & 0 & 0 
    \end{pmatrix},
\end{equation}
the hessian $H_{y}$,
\begin{equation}
    H_{y} = 
    \frac{1}{2} 
    \begin{pmatrix} 
        M_{\tilde{x}xy} \, \partial_{\tilde{x}} x  & (\partial_x m_y^{\uparrow} -\partial_x m_y^{\downarrow}) \partial_{\tilde{x}} x& 0 \\
        (\partial_x m_y^{\uparrow} -\partial_x m_y^{\downarrow}) \partial_{\tilde{x}} x & 0 & 0 \\
        0 & 0 & 0 \\
    \end{pmatrix} 
    =
    \begin{pmatrix} 
        \partial_{\tilde{x}\tilde{x}}y  & \partial_{\tilde{x}\tilde{y}}y & 0 \\
        \partial_{\tilde{y}\tilde{x}}y  & 0 & 0 \\
        0 & 0 & 0 \\
    \end{pmatrix}.
\end{equation}
depends on the entry $\partial_{\tilde{x}} x $ of the nested mapping's Jacobian matrix
and on first partial derivatives of its height functions. In addition, its entry $M_{\tilde{x}xy}$ requires the
second partial derivatives $\partial_{xx} m_y^{\uparrow}$ and $\partial_{xx} m_y^{\downarrow}$, 
\begin{equation}
    M_{\tilde{x}xy} = \partial_{\tilde{x}} M_{xy} = \left((\partial_{xx} m_y^{\uparrow} -\partial_{xx} m_y^{\downarrow}) \tilde{y} 
    + \partial_{xx} m_y^{\uparrow} +\partial_{xx} m_y^{\downarrow} \right) \partial_{\tilde{x}} x
    .
\end{equation}
Finally, the hessian $H_{z} $ consists of entries of $H_y$ and $D M$,
\begin{equation}
    H_{z} = 
    \frac{1}{2}
    \begin{pmatrix} 
        M_{\tilde{x}xz} \partial_{\tilde{x}} x + M_{\tilde{x}yz} \partial_{\tilde{x}} y + M_{yz} \partial_{\tilde{x}\tilde{x}} y & 
            M_{\tilde{x}yz} \partial_{\tilde{y}} y + M_{yz} \partial_{\tilde{x}\tilde{y}} y &   
            M_{\tilde{x}\tilde{z}z} \\
        M_{\tilde{x}yz} \partial_{\tilde{y}} y 
            + M_{yz} \partial_{\tilde{x}\tilde{y}} y &
            M_{\tilde{y}yz} \partial_{\tilde{y}} y & 
            (\partial_y m_z^{\uparrow} -\partial_y m_z^{\downarrow}) \partial_{\tilde{y}} y \\
        M_{\tilde{x}\tilde{z}z} & (\partial_y m_z^{\uparrow} -\partial_y m_z^{\downarrow}) \partial_{\tilde{y}} y & 0 \\
    \end{pmatrix}.
\end{equation}
Its entries $M_{\tilde{x}xz},M_{\tilde{x}yz}, M_{\tilde{y}yz}$ and $M_{\tilde{x}\tilde{z}z}$
\begin{align}
    M_{\tilde{x}xz} &= \partial_{\tilde{x}} M_{xz}(x,y,\tilde{z}) \notag \\
    & = \left((\partial_{xx} m_z^{\uparrow} - \partial_{xx} m_z^{\downarrow}) \tilde{z} 
        + \partial_{xx} m_z^{\uparrow} + \partial_{xx} m_z^{\downarrow} \right) \partial_{\tilde{x}} x
        + \left((\partial_{xy} m_z^{\uparrow} - \partial_{xy} m_z^{\downarrow}) \tilde{z}
        + \partial_{xy} m_z^{\uparrow} + \partial_{xy} m_z^{\downarrow}\right) \partial_{\tilde{x}} y
        ,\\
    M_{\tilde{x}yz} & = \partial_{\tilde{x}} M_{yz}(x,y,\tilde{z}) \notag \\
    & = \left(( \partial_{xy} m_z^{\uparrow} - \partial_{xy} m_z^{\downarrow}) \tilde{z} 
        + \partial_{xy} m_z^{\uparrow} + \partial_{xy} m_z^{\downarrow} \right) \partial_{\tilde{x}}x
        + \left(( \partial_{yy} m_z^{\uparrow} - \partial_{yy} m_z^{\downarrow}) \tilde{z} 
        + \partial_{yy} m_z^{\uparrow} + \partial_{yy} m_z^{\downarrow} \right) \partial_{\tilde{x}}y
        ,\\
    M_{\tilde{y}yz} & = \partial_{\tilde{y}} M_{yz}(x,y,\tilde{z}) \notag \\
    & = \left(( \partial_{xy} m_z^{\uparrow} - \partial_{xy} m_z^{\downarrow}) \tilde{z} 
        + \partial_{xy} m_z^{\uparrow} + \partial_{xy} m_z^{\downarrow} \right) \partial_{\tilde{y}}x
        + \left(( \partial_{yy} m_z^{\uparrow} - \partial_{yy} m_z^{\downarrow}) \tilde{z} 
        + \partial_{yy} m_z^{\uparrow} + \partial_{yy} m_z^{\downarrow} \right) \partial_{\tilde{y}}y
        ,\\
    M_{\tilde{x}\tilde{z}z} & = \partial_{\tilde{x}\tilde{z}} z(x,y,\tilde{z}) \notag \\
    & = (\partial_x m_z^{\uparrow} -\partial_x m_z^{\downarrow}) \partial_{\tilde{x}} x + (\partial_y m_z^{\uparrow} -\partial_y m_z^{\downarrow}) \partial_{\tilde{x}} y 
        ,  
\end{align}
contain the second partial derivatives,
$\partial_{xx} m_z^{\uparrow},$ $\partial_{xx} m_z^{\downarrow}$, 
$\partial_{xy} m_z^{\uparrow} $ and $\partial_{xy} m_z^{\downarrow}$,
of the height functions.

The second derivatives needed in the hessians of $A$ 
are obtained after once more differentiating and rearranging
Equation~\ref{eq:phi_gamma} and
Equation~\ref{eq:phi_dgamma}: 
\begin{align}
    \partial_{xx} a_y & = - \frac{\partial_{xx} \alpha
        + 2 \partial_{xy} \alpha \, \partial_x a_y 
        + \partial_{yy} \alpha \, ( \partial_x a_y)^2 
        }{\partial_y \alpha} 
        ,\quad 
    \partial_{xx} a_z = - \frac{\partial_{xx} \alpha 
        + 2 \partial_{xz} \alpha \,\partial_x a_z 
        + \partial_{zz} \alpha \,( \partial_x a_z)^2 
        }{\partial_z \alpha} 
        ,\\
    \partial_{yy} a_z & = - \frac{\partial_{yy} \alpha 
        + 2 \partial_{yz} \alpha \, \partial_y a_z 
        + \partial_{zz} \alpha \, ( \partial_y a_z)^2 
        }{\partial_z \alpha} 
        ,\quad
    \partial_{xy} a_z  = - \frac{\partial_{xy} \alpha 
        + \partial_{xz} \alpha \, \partial_y a_z 
        + \partial_{yz} \alpha \, \partial_x a_z 
        + \partial_{zz} \alpha \, \partial_x a_z \partial_y a_z  
        }{\partial_z \alpha}.
\end{align}

\end{document}